\documentclass[11pt,reqno]{amsart}

\title[]{Long Time Dynamics of the Three-Dimensional Nernst-Planck-Darcy Model}

\author{Elie Abdo}
\address{Department of Mathematics, University of California Santa Barbara, CA 93106, USA.}
\email{elieabdo@ucsb.edu}

\author{\DJ or\dj e Nikoli\'c}
\address{Department of Mathematics, University of California Santa Barbara, CA 93106, USA.}
\email{nikolic@ucsb.edu}

\usepackage[margin=1in]{geometry}
\usepackage{amsmath, amsthm, amssymb}
\usepackage{times}
\usepackage{comment}
\usepackage{color}
\usepackage{hyperref}
\newcommand{\pa}{\partial}

\newcommand{\la}{\label}
\newcommand{\fr}{\frac}
\newcommand{\na}{\nabla}
\newcommand{\be}{\begin{equation}}
\newcommand{\ee}{\end{equation}}
\newcommand{\ba}{\begin{array}{l}}
\newcommand{\ea}{\end{array}}

\newtheorem{thm}{Theorem}[section]
\newtheorem{prop}{Proposition}[section]

\newcommand{\beg}{\begin}

\newtheorem{lem}{Lemma}[section]
\newtheorem{Thm}{Theorem}[section]

\usepackage{MnSymbol,wasysym}

\newcommand{\R}{\mathbb R}

\def\RR{{\mathbb R}}
\def\TT{{\mathbb T}}
\def\NN{\mathbb N}
\def\PP{\mathbb P}

\begin{document}
\begin{abstract}
We consider an electrodiffusion model describing the evolution of $N$ ionic species in a three-dimensional fluid flowing through a porous medium and forced by added body charges. We address the global well-posedness 
and long-time dynamics of the model. In the absence of added charges, we prove that the ionic concentrations decay exponentially fast in time to their initial spatial averages in all Sobolev norms. When the fluid undergoes the influence of given time-independent charges, we obtain the existence of a finite-dimensional global attractor. 
\end{abstract} 


\maketitle

\tableofcontents

\section{Introduction}

Electrodiffusion in fluids is a 
phenomenon describing the transport of ions through advection, diffusion, and electromigration. It has been extensively studied in the literature in connection with real-world applications, including but not limited to neuroscience \cite{cole1965electrodiffusion, jasielec2021electrodiffusion,koch2004biophysics, lopreore2008computational, mori2009numerical, nicholson2000diffusion, pods2013electrodiffusion,  savtchenko2017electrodiffusion, qian1989electro}, semiconductors 
\cite{biler1994debye,gajewski1986basic,mock1983analysis}, water purification, desalination, and ion separations \cite{gao2014high,lee2016membrane, lee2018diffusiophoretic, yang2019review, zhu2020ion}, and ion-selective membranes \cite{davidson2016dynamical,goldman1989electrodiffusion}.

In this paper, we consider an electrodiffusion model describing the nonlinear and nonlocal evolution of $N$ ionic concentrations $c_i$, with diffusivities $D_i$ and valences $z_i$, interacting with an incompressible fluid, flowing through a porous medium with a velocity $u$ and a pressure $p$. The time evolution of each concentration $c_i$ is described by the Nernst-Planck equations
\be \label{concentrationequation}
\pa_t c_i + \na \cdot j_i = \pa_t c_i + u \cdot \na c_i - D_i \Delta c_i - z_i D_i \na \cdot (c_i \na \Phi) = 0,
\ee for $i = 1, \dots, N$, where $j_i = uc_i - D_i \na c_i - z_iD_i c_i \na \Phi$ are the fluxes, 
$\Phi$ is the electric potential generated by the ionic charge density 
\be \label{densityequation}
\rho = \sum\limits_{i=1}^{N} z_ic_i
\ee and some added time-independent body charge $\tilde{\rho}$, and obeying the following Poisson equation 
\be \label{poissoneq} 
-\Delta \Phi = \rho + \tilde{\rho}.
\ee The velocity of the fluid and electric forces driven by the total potential are balanced via Darcy's law 
\be 
u + \na p = -(\rho + \tilde{\rho}) \na \Phi, \hspace{1cm} \na \cdot u = 0,
\ee which, projected onto the space of divergence-free vectors, reduces to 
\be \label{velocityequation}
u = - \PP((\rho + \tilde{\rho}) \na \Phi).
\ee Here the operator $\PP = I + \na (-\Delta)^{-1} \na \cdot$ denotes the Leray-Hodge projection. 
The model, described by equations \eqref{concentrationequation}--\eqref{velocityequation}, and equipped with prescribed initial data $(c_1(0), \dots, c_N(0))$, is called the Nernst-Planck-Darcy system and will be referred to as NPD. We address the global well-posedness and long-time dynamics of the NPD model for $N$ ionic species with equal diffusivities $D_1 = \dots = D_N :=D$ and valences obeying $|z_1| = \dots = |z_N|$ on the three-dimensional torus $\TT^3 = [0,2\pi]^3$ with periodic boundary conditions.

When the fluid is viscous, Darcy's law is replaced by the Navier-Stokes equations, and the resulting model, called the Nernst-Planck-Navier-Stokes system and abbreviated by NPNS, has been a center of interest over the last two decades. Existence, uniqueness, regularity, analyticity, stability and long-time dynamics of solutions have been widely studied in the literature on two-dimensional domains with and without physical boundaries and under different types of boundary conditions (\cite{abdo2021long, abdo2022space, constantin2019nernst, constantin2021nernst, constantin2022existence, rubinstein2000electro, ryham2009existence, rubinstein2008direct, schmuck2009analysis, zaltzman2007electro} and references theiren). 

The existence of global weak solutions to the NPD system was addressed in the 3D case for two ionic species in \cite{herz2016globalpart1,  herz2012existence} and in the 2D case for $N$ ionic species in \cite{herz2016globalpart2}. In \cite{ignatova2022global}, the authors studied the periodic NPD system for two iconic species with equal diffusivities and respective valences $z$ and $-z$ in the absence of body charges. They proved the existence and uniqueness of global solutions whose spatial Sobolev $H^3$ norms are uniformly bounded in time and whose $L^p$ norms decay exponentially fast time to their initial spatial averages. As a matter of fact, the two-species model has a nice structure by which the nonlinear aspects of the electromigration terms pose fewer major challenges to the analysis of the model. More precisely, the behavior of the ionic concentrations $c_1$ and $c_2$ associated with the two given species is fully determined by that of the charge density $\rho = c_1 - c_2$ and the total amount $\sigma = c_1 + c_2$, in which case the coupled migration 
\be 
\int_{\TT^3} \na \cdot (\sigma \na \Phi) \rho dx + \int_{\TT^3} \na \cdot (\rho \na \Phi) \sigma dx
\ee fully dissipates into $-\|\rho \sqrt{\sigma}\|_{L^2}^2$, yielding good control of the $L^2$ norms of $\rho$ and $\sigma$ and consequently of $c_1$ and $c_2$. However, the case of $N$ ionic species is more challenging and technical. 

In this paper, we consider the three-dimensional periodic NPD system in the more general setting of $N$ ions with equal diffusivities $D$ and valences $|z_1| = \dots = |z_N|$, and we address the long-time dynamics in the presence and absence of body charges $\tilde{\rho}$. The main challenges arise from:
\begin{enumerate}
    \item The more singular constitutive relation between the velocity $u$ and the forces $-\rho \na \Phi$, compared to the viscous Nernst-Planck-Navier-Stokes system;
    \item The nonlinear structure of the electromigration term $\na \cdot (c_i \na \Phi)$ on three-dimensional spatial domains;
    \item The effects of the added time-independent density $\tilde{\rho}$ on the time asymptotic behavior of solutions. 
\end{enumerate}
We successfully overcome these challenges and prove the following results:
\begin{enumerate}
\item Existence of global $C^{\infty}$ solutions that decay exponentially fast in time to their initial spatial averages in all Sobolev norms in the absence of body charges. This generalizes the convergence results obtained in \cite{ignatova2022global} for the $L^p$ norms of two species into the $H^k$ norms of $N$ species for any $k \ge 0$ based on a new bootstrapping argument.    
\item Existence of a smooth finite-dimensional global attractor in the presence of added body charges. 
\end{enumerate}

The paper is organized as follows. In section \ref{s2}, we present some preliminary propositions, including Gronwall lemmas and Poincar\'e inequalities. In section \ref{s3}, we address the existence and uniqueness of a local solution to the 3D NPD system (whose proof is provided in Appendix \ref{AP1}) and we show that it extends globally in time based on a bootstrapping argument. In section \ref{s4}, we present a general framework that summarizes the main conditions required to study the long-time dynamics of nonlinear fractionally dissipative parabolic equations and use it to obtain the smoothness and asymptotic behavior of the NPD model. In section \ref{s5}, we define the solution map associated with the NPD system and address its continuity whose proof is presented in Appendix \ref{AP2}. In section \ref{s6}, we provide a general energy-based argument to obtain the backward uniqueness of solutions to nonlinear parabolic systems and use it to prove the injectivity of the solution map. Finally, the existence of a global attractor in the presence of body charges is shown in Section \ref{s7}. The finite dimensionality of the attractor is based on the decay of volume elements which is established in Appendix \ref{AP3}.

\section{Preliminaries} \la{s2}

For $1 \le p \le \infty$, we denote by $L^p(\TT^3)$ the Lebesgue spaces of periodic measurable functions $f$ from $\TT^3$ to $\R$  such that 
\be 
\|f\|_{L^p} = \left(\int_{\TT^3} |f|^p dx\right)^{1/p} <\infty
\ee if $p \in [1, \infty)$ and 
\be 
\|f\|_{L^{\infty}} = {\mathrm{esssup}}_{\TT^3}  |f| < \infty
\ee if $p = \infty$. We denote the $L^2(\TT^3)$ inner product by $(\cdot,\cdot)_{L^2}$. 

For $k \in \NN$, we denote by $H^k(\TT^3)$ the classical  Sobolev spaces of periodic measurable functions $f$ from $\TT^3$ to $\R$ with weak derivatives of order $k$ such that  
\be 
\|f\|_{H^k}^2 = \sum\limits_{|\alpha| \le k} \|D^{\alpha}f\|_{L^2}^2 < \infty.
\ee 
For a Banach space $(X, \|\cdot\|_{X})$ and $p\in [1,\infty]$, we consider the Lebesgue spaces $ L^p(0,T; X)$ of functions $f$  from $X$ to $\R$  satisfying 
\be 
\int_{0}^{T} \|f\|_{X}^p dt  <\infty,
\ee with the usual convention when $p = \infty$.

We recall the following uniform Gronwall Lemma \cite{abdo2021long}:

\begin{lem} \la{exp} Let $y(t) \geq 0$ obey a differential inequality $$\frac{d}{dt}y + c_1y \leq F_1 + F(t)$$ with initial datum $y(0) = y_0$ with $F_1$ a nonnegative constant, and $F(t) \geq 0$ obeying $$\int\limits_{t}^{t+1} F(s) ds \leq g_0e^{-c_2t} + F_2$$ where $c_1, c_2, g_0$ are positive constants and $F_2$ is a nonnegative constant. Then $$y(t) \leq y_0e^{-c_1t} + g_0e^{c_1+c}(t+1)e^{-ct} + \frac{1}{c_1}F_1 + \frac{e^{c_1}}{1-e^{-c_1}} F_2$$ holds with $c = \min \left\{c_1, c_2 \right\}$.
\end{lem}

The following elementary proposition will be frequently used in the sequel:

\beg{prop} \label{localinteg}
Let $E_1(t), E_2(t), D_1(t)$ and $D_2(t)$ be nonnegative functions defined on $[0, \infty)$, and let $L_1(z)$ and $L_2(z)$ be nonnegative functions defined on $[0,\infty)$
and obeying $L_1(0) = L_2(0) = 0$. For a fixed $z \ge 0$, define $\mathcal{F}_z(t): [0,\infty) \rightarrow [0,\infty)$ by 
\be 
\beg{aligned}
\mathcal{F}_z(t) &= \left[C_1 E_1(t)^{\beta_1} + C_2E_1(t)^{\beta_2} D_1(t)^{\beta_3}+ C_3 L_1(z) + C_4  \right] \\&\quad\quad\quad\quad\quad\quad\times\left[C_5E_2(t)^{\beta_4} + C_6E_2(t)^{\beta_5}D_2(t)^{\beta_6} + C_7L_2(z) \right]
\end{aligned}
\ee where $C_1, \dots, C_7$ are nonnegative universal constants and $\beta_1, \dots, \beta_6$ are nonnegative powers with $\beta_3 +  \beta_6 \in [0,1]$.
If there are positive constants $c_1, \dots, c_4, \kappa_1, \dots, \kappa_4$ and nonnegative functions  $a_1(z), \dots, a_4(z)$ such that $a_1(0) = a_2(0) = a_3(0) = a_4(0) = 0$, 
\be 
E_1(t) \le c_1 e^{-\kappa_1t} + a_1,
\ee 
\be 
E_2(t) \le c_2 e^{-\kappa_2t} + a_2,
\ee 
\be 
\int_{t}^{t+1} D_1(s) ds \le c_3e^{-\kappa_3t} + a_3,
\ee and 
\be 
\int_{t}^{t+1} D_2(s) ds \le c_4e^{-\kappa_4t} + a_4,
\ee for any $t \ge 0$, then there is a positive constant $\kappa $ depending only on $\kappa_1, \dots, \kappa_4$, a positive constant $A_1$ depending only on $c_1, \dots, c_4$, $\kappa_1, \dots, \kappa_4$, $a_1, \dots, a_4$, and a nonnegative function $A_2(z)$ depending only on $a_1, \dots, a_4$ that vanishes at $z=0$, such that the following local-in-time estimate
\be 
\int_{t}^{t+1} \mathcal{F}_z(s) ds 
\le A_1 e^{-\kappa t} + A_2
\ee holds for any $t \ge 0$. 
\end{prop}

\begin{proof}
    The proof of Proposition \ref{localinteg} is straightforward and will be omitted.
\end{proof}

We state and prove a modified version of the  Poincar\'e inequality that will be needed to study the long-time dynamics in $L^p$ spaces:

\beg{prop} (Poincar\'e inequality) Let $p \ge 2$. Let $f \in L^p$ such that $f^{\frac{p-2}{2}}\na f \in L^2$. It holds that  
\be \la{poincareineq}
\|f - \bar{f}\|_{L^p}^p \le Cp^2 \||f- \bar{f}|^{\frac{p-2}{2}} \na f\|_{L^2}^2 
+ \frac{2}{|\TT^3|} \|f - \bar{f}\|_{L^{\fr{p}{2}}}^{p}
\ee where $C$ is a constant independent of $p$.
\end{prop}

\beg{proof}
By the Poincar\'e inequality applied to the mean-free function $(f - \bar{f})^{\frac{p}{2}} - \frac{1}{|\TT^3|} \int_{\TT^3}(f - \bar{f})^{\frac{p}{2}} dx$, we have 
\be
\left\|(f - \bar{f})^{\frac{p}{2}} - \frac{1}{|\TT^3|} \int_{\TT^3}(f - \bar{f})^{\frac{p}{2}} dx \right\|_{L^2}^2
\le C\|\na (f-\bar{f})^{\frac{p}{2}}\|_{L^2}^2
=\frac{Cp^2}{4} \left\|(f-\bar{f})^{\frac{p-2}{2}} \na f\right\|_{L^2}^2
.
\ee By the reverse triangle inequality, we have 
\be 
\beg{aligned}
\left\|(f - \bar{f})^{\frac{p}{2}} - \frac{1}{|\TT^3|} \int_{\TT^3}(f - \bar{f})^{\frac{p}{2}} dx \right\|_{L^2}
&\ge \|(f - \bar{f})^{\frac{p}{2}}\|_{L^2} - \left\| \frac{1}{|\TT^3|} \int_{\TT^3}(f - \bar{f})^{\frac{p}{2}} dx \right\|_{L^2}
\\&= \|(f - \bar{f})^{\frac{p}{2}}\|_{L^2} -\frac{1}{|\TT^3|^{\frac{1}{2}}} \left|\int_{\TT^3}(f - \bar{f})^{\frac{p}{2}} dx \right|.
\end{aligned}
\ee Thus,
\be 
\|f - \bar{f}\|_{L^2}^2 
\le \frac{2Cp^2}{4} \left\|(f-\bar{f})^{\frac{p-2}{2}} \na f\right\|_{L^2}^2
+ \frac{2}{|\TT^3|} \|f - \bar{f}\|_{L^{\frac{p}{2}}}^p.
\ee 
\end{proof}

\section{Existence and Uniqueness of Solutions} \la{s3}

In this section, we prove the existence of unique solutions to the NPD system \eqref{concentrationequation}--\eqref{velocityequation} for $H^1$  initial concentrations.

\beg{defi} A vector $(c_1, \dots, c_N)$ is said to be a strong solution of the NPD system \eqref{concentrationequation}--\eqref{velocityequation} on the time interval $[0,T]$ if it solves \eqref{concentrationequation}--\eqref{velocityequation} in the sense of distributions and obeys the regularity criterion
\be  
c_i \in L^{\infty} (0,T; H^1(\TT^3)) \cap L^2(0,T; H^2(\TT^3))
\ee for all $i \in \left\{1, \dots, N\right\}.$
\end{defi}

\begin{Thm} \la{T1} Let $T > 0$ be an arbitrary time. Let $c_i(0) \in H^1$ for all $i \in \left\{1, \dots, N\right\}$ be nonnegative. Assume that $\int_{\TT^3} (\rho_0 + \tilde{\rho})(x) dx = 0$. The NPD system \eqref{concentrationequation}--\eqref{velocityequation} has a unique strong solution $(c_1, \dots, c_N)$ on $[0,T]$. Moreover, there exist nonnegative constants $C_1, \dots, C_N$ depending only on the initial data, and nonnegative constants $\Gamma_1, \dots, \Gamma_N$ depending only on the body charge $\tilde{\rho}$ and the spatial averages of the initial concentrations, and a positive universal constant $c$ such that the following estimates 
\be 
\|\na c_i(t)\|_{L^2} \le C_i e^{-ct} + \Gamma_i
\ee holds for all positive times $t \ge 0$ and for all $i \in \left\{1, \dots, N\right\}$. Moreover, $\Gamma_1 = \dots = \Gamma_N = 0$ when $\tilde{\rho} = 0$.
\end{Thm}

The proof of Theorem \ref{T1} follows from the following two propositions:

\beg{prop} \la{locso} (Local Solution)
Let $c_i(0) \in H^1$ for all $i \in \left\{1, \dots, N\right\}$.  Assume that $\int_{\TT^3} (\rho_0 + \tilde{\rho})(x) dx = 0$. Then there exists a time $T_0>0$ depending only on the size of the initial data in $H^1$ and the parameters of the problem such that the NPD system \eqref{concentrationequation}--\eqref{velocityequation} has a unique strong solution $(c_1, \dots, c_N)$ on $[0,T_0]$.
\end{prop}

The proof of Proposition \ref{locso} is standard and based on Galerkin approximations, energy estimates, and passage to the limit using the Aubin-Lions lemma. We present the proof in Appendix \ref{AP1} for the sake of completion. 

\beg{prop} (Blow-up Criterion) Let $T > 0$ be an arbitrary time. Let $c_i(0) \in H^1$ for all $i \in \left\{1, \dots, N\right\}$ be nonnegative. Assume that $\int_{\TT^3} (\rho_0 + \tilde{\rho})(x) dx = 0$. Suppose $(c_1, \dots, c_N)$ is a strong solution to the NPD system \eqref{concentrationequation}--\eqref{velocityequation} on $[0,T]$. Then there exist nonnegative constants $C_1, \dots, C_N$ depending only on the initial data, and nonnegative constant $\Gamma_1, \dots, \Gamma_N$ depending only on the body charge $\tilde{\rho}$ and the spatial averages of the initial concentrations, and a positive universal constant $c$ such that the following estimates 
\be 
\|\na c_i(t)\|_{L^2} \le C_i e^{-ct} + \Gamma_i
\ee holds for all positive times $t \ge 0$ and for all $i \in \left\{1, \dots, N\right\}$. Moreover, $\Gamma_1 = \dots = \Gamma_N = 0$ when $\tilde{\rho} = 0$.
\end{prop}

\begin{proof} The nonnegativity of the ionic concentrations is preserved for all times $t \ge 0$, a fact that follows from \cite{constantin2019nernst}.  The proof is divided into several steps.

{\bf{Step 1. Conservation Laws.}} Integrating the ionic concentration equation \eqref{concentrationequation} over $\mathbb{T}^{3}$ yields
\begin{equation}
\frac{d}{dt} \int_{\mathbb{T}^{3}} c_{i}(x,t) dx = 0,
\end{equation} after making use of the divergence-free condition obeyed by the velocity field $u$ and integrating by parts. As a consequence, we have  \begin{equation}
 \int_{\mathbb{T}^{3}} c_{i}(x,t) dx =  \int_{\mathbb{T}^{3}} c_{i}(x,0) dx
\end{equation} for all times $t \ge 0$ and all $i \in \left\{1, \dots, N\right\}$. In view of the relation \eqref{densityequation} between the charge density $\rho$ and the ionic concentrations, we deduce the following conservation law 
 \begin{equation}
 \int_{\mathbb{T}^{3}} \rho (x,t) dx =  \int_{\mathbb{T}^{3}} \rho (x,0) dx
\end{equation} at all times $t \ge 0$.
As the spatial averages of the charge density and ionic concentrations are time-independent, we respectively denote them \begin{align}
 \bar{\rho} &= \frac{1}{|\mathbb{T}^{3}|} \int_{\mathbb{T}^{3}} \rho(x,t) dx, \\
\bar{c_{i}} &= \frac{1}{|\mathbb{T}^{3}|} \int_{\mathbb{T}^{3}} c_{i}(x,t) dx,
\end{align} for all $i \in \left\{1, \dots, N \right\}$.

{\bf{Step 2. A priori Bounds for the Charge Density in $L^2$.}} Let $\rho = \sum\limits_{i=1}^{N}z_{i}c_{i}$ and $\sigma = \sum\limits_{i=1}^{N}c_{i}.$ 
Then the time evolutions of $\rho$ and $\sigma$ are described by \begin{equation} \label{something}
    \partial_{t}\rho + u \cdot \nabla \rho - D\Delta \rho = Dz^{2}\nabla \cdot (\sigma \nabla \Phi)
\end{equation} and 
\be 
\label{something2}
    \partial_{t}(z^{2}\sigma) + u \cdot \nabla(z^{2}\sigma) - D\Delta(z^{2}\sigma)= Dz^{2}\nabla \cdot (\rho\nabla \Phi),
\ee respectively. We multiply the  $\rho$-equation \eqref{something} by $\rho$  and the $\sigma$-equation $\eqref{something2}$ by $\sigma$, integrate spatially over $\mathbb{T}^{3}$, and obtain the energy equalities
\be \label{prop321}
\frac{1}{2} \frac{d}{dt} \|\rho\|_{L^{2}}^{2} + D \|\nabla \rho\|_{L^{2}}^{2} = -Dz^{2} \int_{\mathbb{T}^{3}} \sigma \nabla \Phi \cdot \nabla \rho dx
\ee  and
\be \label{prop322}
\frac{z^{2}}{2} \frac{d}{dt} \|\sigma\|_{L^{2}}^{2} + Dz^2 \|\nabla \sigma\|_{L^{2}}^{2} = Dz^{2} \int_{\mathbb{T}^{3}}  (\nabla \rho \cdot \nabla \Phi) \sigma dx + Dz^{2} \int_{\mathbb{T}^{3}} \rho \Delta \Phi \sigma dx
\ee respectively. We add \eqref{prop321} and \eqref{prop322}, make use of the Poisson equation \eqref{poissoneq} obeyed by $\Phi$, and deduce the differential equation 
\be 
\fr{1}{2} \fr{d}{dt} \left(\|\rho\|_{L^2}^2 + z^2 \|\sigma\|_{L^2}^2 \right)  + D\left(\|\na \rho\|_{L^2}^2 + z^2 \|\na \sigma\|_{L^2}^2 \right)
= -Dz^2 \int_{\TT^3} \rho (\rho + \tilde{\rho}) \sigma dx,
\ee which boils down to 
\be \label{prop325}
\fr{1}{2} \fr{d}{dt} \left(\|\rho\|_{L^2}^2 + z^2 \|\sigma\|_{L^2}^2 \right)  + D\left(\|\na \rho\|_{L^2}^2 + z^2 \|\na \sigma\|_{L^2}^2 + z^2 \|\rho \sqrt{\sigma}\|_{L^2}^2  \right) = - Dz^2 \int_{\TT^3} \rho \tilde{\rho} \sigma dx
\ee due to the nonnegativity of the ionic concentrations.  
An application of H\"older and Young inequalities gives 
\be 
\left|Dz^2 \int_{\TT^3} \rho \tilde{\rho} \sigma dx \right|
\le Dz^2 \int_{\TT^3} |\rho \sqrt{\sigma}| |\rho \sqrt{\sigma}|dx 
\le \frac{Dz^2}{2} \|\rho \sqrt{\sigma}\|_{L^2}^2 + \frac{Dz^2}{2} \|\tilde{\rho} \sqrt{\sigma}\|_{L^2}^2.
\ee   Furthermore, we estimate \be 
\|\tilde{\rho} \sqrt{\sigma}\|_{L^2}^2
= \int_{\TT^3} \tilde{\rho}^2 (\sigma - \bar{\sigma}) dx 
+ \int_{\TT^3} \tilde{\rho}^2 \bar{\sigma} dx 
\le \|\tilde{\rho}^2\|_{L^2} \|\sigma - \bar{\sigma}\|_{L^2} +\bar{\sigma} \|\tilde{\rho}\|_{L^2}^2,
\ee where $\bar{\sigma}$ stands for the time-independent spatial average of $\sigma$ over $\TT^3$. In view of the Poincar\'e inequality applied to the mean-free function $\sigma - \bar{\sigma}$ and Young's inequality for products, we have 
\be \la{prop326}
\|\tilde{\rho} \sqrt{\sigma}\|_{L^2}^2
\le C\|\tilde{\rho}^2\|_{L^2} \|\na \sigma\|_{L^2} + \bar{\sigma} \|\tilde{\rho}\|_{L^2}^2
\le \frac{Dz^2}{2} \|\na \sigma\|_{L^2}^2 + C \|\tilde{\rho}\|_{L^4}^4 + \bar{\sigma} \|\tilde{\rho}\|_{L^2}^2.
\ee Putting \eqref{prop325}--\eqref{prop326} together, we infer that 
\be 
\fr{d}{dt} \left(\|\rho\|_{L^2}^2 + z^2 \|\sigma\|_{L^2}^2 \right) + D \left(\|\na \rho\|_{L^2}^2 + z^2 \|\na \sigma\|_{L^2}^2 + z^2 \|\rho \sqrt{\sigma}\|_{L^2}^2 \right)
\le C\|\tilde{\rho}\|_{L^4}^4 + \bar{\sigma} \|\tilde{\rho}\|_{L^2}^2.
\ee As $\bar{\rho}$ and $\bar{\sigma}$ are constants in time, it holds that 
\be 
\fr{d}{dt} \|\rho\|_{L^2}^2 = \fr{d}{dt} \|\rho - \bar{\rho}\|_{L^2}^2 
\ee and 
\be 
\fr{d}{dt} \|\sigma\|_{L^2}^2 = \frac{d}{dt} \|\sigma - \bar{\sigma}\|_{L^2}^2,
\ee from which we obtain the differential inequality
\be \label{prop328}
\fr{d}{dt} \left(\|\rho - \bar{\rho}\|_{L^2}^2 + z^2 \|\sigma - \bar{\sigma}\|_{L^2}^2 \right) + D \left(\|\na \rho\|_{L^2}^2 + z^2 \|\na \sigma\|_{L^2}^2 + z^2 \|\rho \sqrt{\sigma}\|_{L^2}^2 \right)
\le C\|\tilde{\rho}\|_{L^4}^4 + \bar{\sigma} \|\tilde{\rho}\|_{L^2}^2.
\ee By the Poincar\'e inequality, it follows that 
\be 
\fr{d}{dt} \left(\|\rho - \bar{\rho}\|_{L^2}^2 + z^2 \|\sigma - \bar{\sigma}\|_{L^2}^2 \right) + cD \left(\|\rho - \bar{\rho}\|_{L^2}^2 + z^2 \|\sigma - \bar{\sigma}\|_{L^2}^2  \right)
\le C\|\tilde{\rho}\|_{L^4}^4 + \bar{\sigma} \|\tilde{\rho}\|_{L^2}^2.
\ee Finally, we apply the uniform Gronwall Lemma \ref{exp} and deduce that 
\be \la{prop327}
\|\rho(t) - \bar{\rho}\|_{L^2}^2 + z^2 \|\sigma(t)- \bar{\sigma}\|_{L^2}^2 
\le A_1 e^{-cDt} + A_2
\ee holds for all $t \ge 0$, where $A_1 = \|\rho_0 - \bar{\rho}\|_{L^2}^2 + z^2 \|\sigma_0 - \bar{\sigma}\|_{L^2}^2$ and $A_2 = C(z,D) \left(\|\tilde{\rho}\|_{L^4}^4 + \|\tilde{\rho}\|_{L^2}^2 \right)$.
We point out that in the case of two ionic species ($N=2$) with valences $z_1 = z>0$ and $z_2 = -z$, we have the instantaneous estimates
\be 
\|c_1(t) - \bar{c}_1\|_{L^2}^2 + \|c_2(t) - \bar{c}_2\|_{L^2}^2 
\le A_1e^{-cDt} + A_2
\ee due to \eqref{prop327} and the parallelogram law. In the more general case of $N$ ionic species, we need to study the evolution of each ionic concentration. 

{\bf{Step 3. A priori Bounds for the Concentrations in $L^2$.}} We take the $L^2$ inner product of the ionic concentration equation \eqref{concentrationequation} with $c_i$, make use of the divergence-free condition obeyed by $u$, integrate by parts, and obtain the energy equality
\be 
\frac{1}{2} \frac{d}{dt} \| c_{i}\|_{L^{2}}^{2} + D \| \nabla c_{i} \|_{L^{2}}^{2} = -Dz_{i} \int_{\mathbb{T}^{3}} c_{i}\nabla \Phi \nabla c_{i}dx.
\ee In view of the nonnegativity of the ionic concentrations, we estimate 
\be 
\beg{aligned}
 &\left|\int_{\mathbb{T}^{3}} Dz_{i}c_{i}\nabla \Phi \nabla c_{i}dx \right| 
 \le D|z_i| \int_{\TT^3} c_i |\na \Phi| |\na c_i|  dx
 \le D|z_i| \int_{\TT^3} \left( \sum\limits_{j=1}^{N} c_j \right) |\na \Phi| |\na c_i|
 \\&=  D|z_i| \int_{\TT^3} \sigma  |\na \Phi| |\na c_i| dx
 =  D|z_i| \int_{\TT^3} \left(\sigma - \bar{\sigma}\right)  |\na \Phi| |\na c_i| dx
 + D|z_i| \bar{\sigma} \int_{\TT^3} |\na \Phi| |\na c_i| dx,
 \end{aligned}
 \ee which reduces to
\be 
 \left|\int_{\mathbb{T}^{3}} Dz_{i}c_{i}\nabla \Phi \nabla c_{i}dx \right| \leq \frac{D}{2} \| \nabla c_{i} \|_{L^{2}}^{2} +D|z_{i}|^{2} \| \nabla \Phi \|_{L^{\infty}}^{2}\| \sigma - \bar{\sigma} \|_{L^{2}}^{2} + D|z_{i}|^{2} \bar{\sigma}^{2} \| \nabla \Phi \|_{L^{2}}^{2}  
\ee after straightforward applications of H\"older and Young inequalities. Elliptic regularity estimates yield
\be 
\| \nabla \Phi \|_{L^{\infty}} \leq C \| \nabla \rho + \nabla \tilde{\rho} \|_{L^{2}} \leq C \| \nabla \rho \|_{L^{2}} +  C \| \nabla \tilde{\rho}\|_{L^{2}}.
\ee Consequently, the time rate of change of the $L^2$ norm of each concentration satisfies the energy inequality 
\be 
\frac{1}{2} \frac{d}{dt} \| c_{i} - \bar{c}_i\|_{L^{2}}^{2} + \frac{D}{2} \| \nabla c_{i}\|_{L^{2}}^{2}\leq CD|z|^{2} (1+\bar{\sigma}^{2}) ( \|\na \rho \|_{L^{2}}^2 +  \|\na \tilde{\rho} \|_{L^{2}}^2) .  
\ee Integrating \eqref{prop328} in time from $t$ to $t+1$, we have 
\be 
D \int_{t}^{t+1} \|\na \rho(s)\|_{L^2}^2 ds 
\le A_1 e^{-cDt} + A_2 + C\|\tilde{\rho}\|_{L^4}^4 + \bar{\sigma} \|\tilde{\rho}\|_{L^2}^2  
\ee for any $t \ge 0$. As the conditions of the uniform Gronwall Lemma \ref{exp} are satisfied, we deduce the existence of a nonnegative constant $A_3$ depending only on the $L^2$ norm of the initial concentrations, a nonnegative constant $A_4$ depending only on the $H^1$ norm of $\tilde{\rho}$ (with $A_4 =0$ if $\tilde{\rho} = 0$, and a positive universal constant $c_1$ such that 
\be 
\|c_i(t) - \bar{c}_i\|_{L^2}^2 
\le A_3 e^{-c_1Dt} + A_4
\ee  holds for all $t \ge 0$ and $i \in \left\{1, \dots, N \right\}$.

{\bf{Step 4. A priori Bounds for the Concentrations in $L^p$ spaces for $p \in \left\{3, 4, 6\right\}$.}}  The densities $\rho - \bar{\rho}$ and $\sigma - \bar{\sigma}$ evolve according to 
\be \la{step411}
\partial_{t} (\rho - \bar{\rho}) + u\na(\rho - \bar{\rho}) - D\Delta(\rho - \bar{\rho}) = Dz^{2} \na \cdot (\sigma \na \phi)
\ee and 
\be \la{step412}
\partial_{t} (\sigma - \bar{\sigma}) + u\na(\sigma - \bar{\sigma}) - D\Delta(\sigma - \bar{\sigma}) = D \na \cdot (\rho \na \phi)
\ee respectively. We multiply the $\rho$-equation \eqref{step411} by $(\rho - \bar{\rho})|\rho - \bar{\rho}|^{p-2}$ and the $\sigma$-equation \eqref{step412} by $(\sigma - \bar{\sigma})|\sigma - \bar{\sigma}|^{p-2}$ and integrate over $\TT^3$. We obtain the following energy equalities
\be 
\fr{1}{p}\fr{d}{dt}\|\rho - \bar{\rho}\|_{L^{p}}^{p} + D(p-1) \| |\rho - \bar{\rho}|^{\fr{p-2}{2}}\na \rho \|_{L^{2}}^{2} = Dz^2  \int_{\TT^3} \na \cdot (\sigma \na \Phi) (\rho - \bar{\rho}) |\rho - \bar{\rho}|^{p-2} dx
\ee and 
\be 
\fr{1}{p}\fr{d}{dt}\|\sigma - \bar{\sigma}\|_{L^{p}}^{p} + D(p-1) \| |\sigma - \bar{\sigma}|^{\fr{p-2}{2}}\na \sigma \|_{L^{2}}^{2} = D  \int_{\TT^3} \na \cdot (\rho \na \Phi) (\sigma - \bar{\sigma}) |\sigma - \bar{\sigma}|^{p-2} dx. 
\ee An application of H\"older and Young inequalities  give 
\be 
\beg{aligned}
&\left|Dz^2  \int_{\TT^3} \na \cdot (\sigma \na \Phi) (\rho - \bar{\rho}) |\rho - \bar{\rho}|^{p-2} dx \right| 
= \left|Dz^2 (p-1) \int_{\TT^3}  (\sigma \na \Phi) \cdot \na \rho  |\rho - \bar{\rho}|^{p-2} dx \ \right|
\\&\quad\quad\le Dz^{2}(p-1)\| (\na \rho)(\rho - \bar{\rho})^{\fr{p-2}{2}} \|_{L^{2}} \|\sigma \na \phi\|_{L^{p}}\|\rho - \bar{\rho}\|_{L^{p}}^{\fr{p-2}{2}} \\&\quad\quad\le \fr{D(p-1)}{2} \|\na \rho (\rho - \bar{\rho})^{\fr{p-2}{2}} \|_{L^{2}}^{2} + C_{z,p,D} \|\sigma \na \Phi \|_{L^{p}}^{2}\|\rho - \bar{\rho}\|_{L^{p}}^{p-2}
\end{aligned}
\ee and 
\be 
\beg{aligned}
&\left|D  \int_{\TT^3} \na \cdot (\rho \na \Phi) (\sigma - \bar{\sigma}) |\sigma - \bar{\sigma}|^{p-2} dx \right| 
= \left|D (p-1) \int_{\TT^3}  (\rho \na \Phi) \cdot \na \sigma |\sigma - \bar{\sigma}|^{p-2} dx \ \right|
\\&\quad\quad\le D(p-1)\| (\na \sigma)(\sigma - \bar{\sigma})^{\fr{p-2}{2}} \|_{L^{2}} \|\rho \na \Phi\|_{L^{p}} \|\sigma - \bar{\sigma}\|_{L^{p}}^{\fr{p-2}{2}} 
\\&\quad\quad\le \fr{D(p-1)}{2} \|\na \sigma (\sigma - \bar{\sigma})^{\fr{p-2}{2}} \|_{L^{2}}^{2} + C_{p,D}  \|\rho \na \Phi\|_{L^{p}}^{2}\|\sigma - \bar{\sigma}\|_{L^{p}}^{p-2}
\end{aligned}
\ee Consequently, we obtain the differential inequalities
\be  \label{4.1.1}
\fr{1}{p}\fr{d}{dt}\|\rho - \bar{\rho}\|_{L^{p}}^{p} + \frac{cD(p-1)}{p^2} \|\rho - \bar{\rho} \|_{L^{p}}^{p} \le C_{z,p,D}  \|\sigma \na \Phi \|_{L^{p}}^{2}\|\rho - \bar{\rho}\|_{L^{p}}^{p-2} + C\|\rho - \bar{\rho}\|_{L^{\fr{p}{2}}}^{p}
\ee and 
\be \label{4.1.2.}
\fr{1}{p}\fr{d}{dt}\|\sigma - \bar{\sigma}\|_{L^{p}}^{p} + \fr{cD(p-1)}{p^2} \|\sigma - \bar{\sigma} \|_{L^{p}}^{p} \le  C_{z,p,D} \|\rho \na \Phi \|_{L^{p}}^{2}\|\sigma - \bar{\sigma}\|_{L^{p}}^{p-2} + C\|\sigma - \bar{\sigma}\|_{L^{\fr{p}{2}}}^{p}
\ee  after applying the Poincar\'e inequality \eqref{poincareineq}.

{\bf{Substep 4.1. A priori Bounds for the Concentrations in $L^3$.}} We set $p=3$ in \eqref{4.1.1} and \eqref{4.1.2.}. We bound 
\be 
\|\rho - \bar{\rho}\|_{L^{\fr{3}{2}}}^{3}
\le C\|\rho - \bar{\rho}\|_{L^2}^{2} \|\rho - \bar{\rho}\|_{L^3}
\ee and 
\be 
\|\sigma - \bar{\sigma}\|_{L^{\fr{3}{2}}}^{3}
\le C\|\sigma - \bar{\sigma}\|_{L^2}^{2} \|\sigma - \bar{\sigma}\|_{L^3}
\ee and obtain the evolution inequalities
\be \la{4.1.5}
\fr{d}{dt} \|\rho - \bar{\rho}\|_{L^3}^2 + \frac{2cD}{9} \|\rho - \bar{\rho}\|_{L^3}^2 
\le C_{z,D} \|\sigma \na \Phi\|_{L^3}^2 + C\|\rho - \bar{\rho}\|_{L^2}^{2}
\ee and 
\be \la{4.1.4}
\fr{d}{dt} \|\sigma - \bar{\sigma}\|_{L^3}^2 + \frac{2cD}{9} \|\sigma - \bar{\sigma}\|_{L^3}^2 
\le C_{z,D} \|\rho \na \Phi\|_{L^3}^2 + C\|\sigma- \bar{\sigma}\|_{L^2}^{2}.
\ee We add \eqref{4.1.5} and \eqref{4.1.4}. Using elliptic estimates and Sobolev embeddings, we have 
\be 
\beg{aligned}
\|\sigma \na \Phi\|_{L^3}^2 + \|\rho \na \Phi\|_{L^3}^2
&\le C\left(\|\sigma\|_{L^6}^2 + \|\rho\|_{L^6}^2\right) \|\na \Phi\|_{L^6}^2
\\&\le C\left(\|\sigma - \bar{\sigma}\|_{L^6}^2 + \|\rho - \bar{\rho}\|_{L^6}^2 + \bar{\sigma}^2 + \bar{\rho}^2 \right) \|\na \Phi\|_{L^6}^2 
\\&\le C\left(\|\na \sigma\|_{L^2}^2 + \|\na \rho\|_{L^2}^2 + \bar{\sigma}^2 + \bar{\rho}^2 \right) \|\rho + \tilde{\rho}\|_{L^2}^2
\\&\le C\left(\|\na \sigma\|_{L^2}^2 + \|\na \rho\|_{L^2}^2 + \bar{\sigma}^2 + \bar{\rho}^2 \right) \left(\|\rho - \bar{\rho}\|_{L^2}^2 + \|\tilde{\rho} - \bar{\rho}\|_{L^2}^2\right),
\end{aligned}
\ee which gives rise to the differential inequality 
\be 
\beg{aligned}
\frac{d}{dt} &\left(\|\rho - \bar{\rho}\|_{L^3}^2 + \|\sigma - \bar{\sigma}\|_{L^3}^2 \right) + \frac{2cD}{9} \left(\|\rho - \bar{\rho}\|_{L^3}^2 + \|\sigma - \bar{\sigma}\|_{L^3}^2 \right)
\\&\quad\quad\le C \left(\|\na \sigma\|_{L^2}^2 + \|\na \rho\|_{L^2}^2 + \bar{\sigma}^2 + \bar{\rho}^2 \right) \left(\|\rho - \bar{\rho}\|_{L^2}^2 + \|\tilde{\rho} - \bar{\rho}\|_{L^2}^2\right) + \|\rho - \bar{\rho}\|_{L^2}^{{2}} + \|\sigma - \bar{\sigma}\|_{L^2}^{{2}}.
\end{aligned}
\ee We note that $\|\tilde{\rho} - \bar{\rho}\|_{L^2}$ vanishes when $\tilde{\rho} =0$. Indeed, integrating the Poisson equation obeyed by $\Phi$ over $\TT^3$ gives
\be 
\int_{\TT^3} (\tilde{\rho} + \rho) dx = -\int_{\TT^3} \Delta \Phi dx = 0. 
\ee As a consequence of Step 3, Proposition \ref{localinteg}, and Lemma \ref{exp},  there are positive constants $A_5$ and $c_2$ and a nonnegative constant $A_6$ such that 
\be 
\|\rho(t) - \bar{\rho}\|_{L^3}^2 + \|\sigma(t) - \bar{\sigma}\|_{L^3}^2
\le A_5 e^{-c_2Dt} + A_6
\ee for all $t \ge 0$, where $A_6$ vanishes when $\tilde{\rho} = 0$.

{\bf{Substep 4.2. A priori Bounds for the Concentrations in $L^4$.}} The $L^4$ norm of $\|\rho - \bar{\rho}\|_{L^4}^2 + \|\sigma - \bar{\sigma}\|_{L^4}^2$ obeys
\be 
\beg{aligned}
\frac{d}{dt} &\left(\|\rho - \bar{\rho}\|_{L^4}^2 +\|\sigma - \bar{\sigma}\|_{L^4}^2  \right)
+ \frac{3cD}{16} \left(\|\rho - \bar{\rho}\|_{L^4}^2 +\|\sigma - \bar{\sigma}\|_{L^4}^2  \right)
\\&\quad\quad\le C_{z,D} \left(\|\sigma \na \Phi\|_{L^4}^2 + \|\rho \na \Phi\|_{L^4}^2 \right) + C\left(\|\rho - \bar{\rho}\|_{L^2}^{{2}} + \|\sigma - \bar{\sigma}\|_{L^2}^{2} \right).
\end{aligned}
\ee The 3D elliptic estimate 
\be 
\|\na \Phi\|_{L^{\infty}}
= \|\na (-\Delta)^{-1} (\rho + \tilde{\rho})\|_{L^{\infty}}
\le C\|\na \na (-\Delta)^{-1} (\rho + \tilde{\rho})\|_{L^4}
\le C\|\rho + \tilde{\rho}\|_{L^4}
\ee followed by 3D Gagliardo-Nirenberg interpolation inequalities
\be 
\|\rho + \tilde{\rho}\|_{L^4}^2
\le C\|\rho - \bar{\rho}\|_{L^4}^2 + C\|\tilde{\rho} -\bar{\rho}\|_{L^4}^2
\le C\|\rho - \bar{\rho}\|_{L^3}\|\na \rho\|_{L^2} + C\|\tilde{\rho}- \bar{\rho}\|_{L^4}^2,
\ee 
\be 
\|\sigma\|_{L^4}^2 
\le C\|\sigma - \bar{\sigma}\|_{L^4}^2 + C\bar{\sigma}^2
\le C\|\sigma - \bar{\sigma}\|_{L^3} \|\na \sigma\|_{L^2} + C\bar{\sigma}^2,
\ee and 
\be 
\|\rho\|_{L^4}^2 
\le C\|\rho - \bar{\rho}\|_{L^4}^2 + C\bar{\rho}^2 
\le C\|\rho - \bar{\rho}\|_{L^3}\|\na \rho\|_{L^2} + C\bar{\rho}^2
\ee yield the bound 
\be 
\beg{aligned}
&\|\sigma \na \Phi\|_{L^4}^2 + \|\rho \na \Phi\|_{L^4}^2
\le C\left(\|\sigma\|_{L^4}^2+ \|\rho\|_{L^4}^2 \right) \|\na \Phi\|_{L^{\infty}}^2
\\&\quad\quad\le C\left(\|\sigma -\bar{\sigma}\|_{L^3}\|\na \sigma\|_{L^2} + \|\rho - \bar{\rho}\|_{L^3} \|\na \rho\|_{L^2} + \bar{\sigma}^2 + \bar{\rho}^2 \right)\left(\|\rho - \bar{\rho}\|_{L^3} \|\na \rho\|_{L^2} + \|\tilde{\rho} - \bar{\rho}\|_{L^4}^2 \right).
\end{aligned}
\ee By making use of Substep 4.1, Proposition \ref{localinteg}, and Lemma \ref{exp}, we infer the existence of positive constants $A_7$ and $c_3$ and a nonnegative constant $A_8$ such that 
\be 
\|\rho(t) - \bar{\rho}\|_{L^4}^2 + \|\sigma(t) - \bar{\sigma}\|_{L^4}^2
\le A_7 e^{-c_3Dt} + A_8
\ee for all $t \ge 0$ where $A_8$ vanishes when $\tilde{\rho} = 0$.

{\bf{Substep 4.3. A priori Bounds for the Concentrations in $L^6$.}} The $L^6$ norm of $\|\rho - \bar{\rho}\|_{L^6}^2 + \|\sigma - \bar{\sigma}\|_{L^6}^2$ obeys
\be \la{4.3.2}
\beg{aligned}
\frac{d}{dt} &\left(\|\rho - \bar{\rho}\|_{L^6}^2 +\|\sigma - \bar{\sigma}\|_{L^6}^2  \right)
+ \frac{5cD}{36} \left(\|\rho - \bar{\rho}\|_{L^6}^2 +\|\sigma - \bar{\sigma}\|_{L^6}^2  \right)
\\&\quad\quad\le C_{z,D} \left(\|\sigma \na \Phi\|_{L^6}^2 + \|\rho \na \Phi\|_{L^6}^2 \right) + C\left(\|\rho - \bar{\rho}\|_{L^3}^{{2}} + \|\sigma - \bar{\sigma}\|_{L^3}^{{2}} \right).
\end{aligned}
\ee Elliptic regularity and continuous Sobolev embeddings give rise to 
\be \la{4.3.1} 
\beg{aligned}
&\|\sigma \na \Phi\|_{L^6}^2 + \|\rho \na \Phi\|_{L^6}^2
\le \left(\|\sigma\|_{L^6}^2 + \|\rho\|_{L^6}^2 \right) \|\na \Phi\|_{L^{\infty}}^2
\\&\quad\quad\le C\left(\|\sigma - \bar{\sigma}\|_{L^6}^2 + \|\rho - \bar{\rho}\|_{L^6}^2 + \bar{\sigma}^2 + \bar{\rho}^2 \right) \|\rho + \tilde{\rho}\|_{L^4}^2
\\&\quad\quad\le C\left(\|\na \sigma\|_{L^2}^2 + \|\na \rho\|_{L^2}^2 + \bar{\sigma}^2 + \bar{\rho}^2 \right)\left(\|\rho - \bar{\rho}\|_{L^4}^2 + \|\tilde{\rho}- \bar{\rho}\|_{L^4}^2 \right).
\end{aligned}
\ee Putting  \eqref{4.3.1} and \eqref{4.3.2} together, using the $L^4$ estimates derived in Step 4.2, and applying Proposition \ref{localinteg} and Lemma \ref{exp}, we deduce the existence of positive constants $A_9$ and $c_4$ and a nonnegative constant $A_{10}$ such that 
\be 
\|\rho(t) - \bar{\rho}\|_{L^6}^2 + \|\sigma(t) - \bar{\sigma}\|_{L^6}^2 
\le A_9 e^{-c_4Dt} + A_{10}
\ee for all $t \ge 0$, where $A_{10}$ vanishes when $\tilde{\rho} = 0$. 

{\bf{Step 5. A priori Bounds for the Concentration Gradients in $L^2$.}} The $L^2$ norm of $\na c_i$ evolves according to 
\be 
\frac{1}{2} \frac{d}{dt} \|\na c_i\|_{L^2}^2
+ D \|\Delta c_i\|_{L^2}^2
= \int_{\TT^3} u \cdot \na c_i \Delta c_i dx
+ z_iD \int_{\TT^3} \na \cdot (c_i \na\Phi) \Delta c_i dx.
\ee Applications of the H\"older, Young, and Gagliardo-Nirenberg inequalities yield the bounds
\be 
\beg{aligned}
\left|\int_{\TT^3} u \cdot \na c_i \Delta c_i dx \right| 
&\le \|u\|_{L^6} \|\na c_i\|_{L^3} \|\Delta c_i\|_{L^2}
\le C\|u\|_{L^6}\|\na c_i\|_{L^2}^{\frac{1}{2}}\|\Delta c_i\|_{L^2}^{\frac{3}{2}}
\\&\le \frac{D}{4} \|\Delta c_i\|_{L^2}^2
+ C\|u\|_{L^6}^4 \|\na c_i\|_{L^2}^2
\end{aligned}
\ee and 
\be 
\beg{aligned}
&D|z_i|\left|\int_{\TT^3} \na \cdot (c_i \na \Phi) \Delta c_i dx \right|
\\&\quad\quad\le D|z_i| \int_{\TT^3} |\na c_i| |\na \Phi| |\Delta c_i| dx 
+ D|z_i| \int_{\TT^3} |c_i - \bar{c}_i| |\Delta \Phi| |\Delta c_i| dx 
+ D|z_i| |\bar{c}_i| \int_{\TT^3} |\Delta \Phi| |\Delta c_i| dx
\\&\quad\quad\le CD|z_i| \|\Delta c_i\|_{L^2} \left(\|\na c_i\|_{L^2} \|\na \Phi\|_{L^{\infty}} + \|c_i - \bar{c}_i\|_{L^6} \|\Delta \Phi\|_{L^3} + |\bar{c}_i| \|\Delta \Phi\|_{L^2} \right)
\\&\quad\quad\le CD|z_i| \|\Delta c_i\|_{L^2} \left(\|\na c_i\|_{L^2} \|\na \Phi\|_{L^{\infty}} + \|\na c_i\|_{L^2} \|\Delta \Phi\|_{L^3} + |\bar{c}_i| \|\Delta \Phi\|_{L^2} \right)
\\&\quad\quad\le \frac{D}{4}\|\Delta c_i\|_{L^2}^2
+ C\|\na c_i\|_{L^2}^2 \left(\|\na \Phi\|_{L^{\infty}}^2 + \|\Delta \Phi\|_{L^3}^2 \right) + C\bar{c}_i^2 \|\Delta \Phi\|_{L^2}^2.
\end{aligned}
\ee Thus, we obtain 
\be 
\beg{aligned}
\frac{d}{dt} \|\na c_i\|_{L^2}^2
+ D\|\Delta c_i\|_{L^2}^2
\le C\|\na c_i\|_{L^2}^2 \left(\|u\|_{L^6}^4 + \|\na \Phi\|_{L^{\infty}}^2 + \|\Delta \Phi\|_{L^3}^2 \right) + C\bar{c}_i^2 \|\Delta \Phi\|_{L^2}^2.
\end{aligned}
\ee In view of the boundedness of the Leray projector $\PP$ on $L^6$, we estimate
\be 
\|u\|_{L^6}^4 
\le \|\PP ((\rho + \tilde{\rho})\na\Phi)\|_{L^6}^4
\le C\|\rho + \tilde{\rho}\|_{L^6}^4 \|\na \Phi\|_{L^{\infty}}^4
\le C\|\rho + \tilde{\rho}\|_{L^6}^8
\le C\|\rho - \bar{\rho}\|_{L^6}^8 + C\|\bar{\rho} - \tilde{\rho}\|_{L^6}^8.
\ee Using the Poisson equation obeyed by $\Phi$, we have 
\be 
\|\na \Phi\|_{L^{\infty}}^2
\le C\|\rho + \tilde{\rho}\|_{L^4}^2
\le C\|\rho - \bar{\rho}\|_{L^4}^2 + C\|\tilde{\rho} - \bar{\rho}\|_{L^4}^2,
\ee 
\be 
\|\Delta \Phi\|_{L^3}^2
\le \|\rho + \tilde{\rho}\|_{L^3}^2
\le C\|\rho - \bar{\rho}\|_{L^3}^2 + C\|\tilde{\rho} - \bar{\rho}\|_{L^3}^2,
\ee and 
\be 
\|\Delta \Phi\|_{L^2}^2
\le \|\rho + \tilde{\rho}\|_{L^2}^2
\le C\|\rho - \bar{\rho}\|_{L^2}^2 + C\|\tilde{\rho} - \bar{\rho}\|_{L^2}^2.
\ee Therefore, the evolution of $\|\na c_i\|_{L^2}$ boils down to 
\be 
\beg{aligned}
\frac{d}{dt} \|\na c_i\|_{L^2}^2 + D\|\Delta c_i\|_{L^2}^2
&\le C\|\na c_i\|_{L^2}^2 \left(\|\rho - \bar{\rho}\|_{L^6}^8 + \|\tilde \rho - \bar{\rho}\|_{L^6}^8 + \|\rho - \bar{\rho}\|_{L^4}^2 + \|\bar{\rho} - \tilde{\rho}\|_{L^4}^2 \right) 
\\&\quad\quad+ C\bar{c}_i^2 \left(\|\rho - \bar{\rho}\|_{L^2}^2 + \|\tilde{\rho} - \bar{\rho}\|_{L^2}^2 \right).
\end{aligned}
\ee By steps 3 and 4, we deduce that the right-hand side of this latter energy inequality obeys the conditions of Proposition \ref{localinteg} and consequently, we can apply the Gronwall Lemma \ref{exp} to obtain positive constants $A_{11}$ and $c_5$ and a nonnegative constant $A_{12}$ such that 
\be 
\sum\limits_{i=1}^{N} \|\na c_i(t)\|_{L^2}^2 \le A_{11} e^{-c_5Dt} + A_{12}
\ee for all $t \ge 0$, where $A_{12}$ vanishes when $\tilde{\rho} =0$. 
    
\end{proof}

\section{Higher Regularity of Solutions} \la{s4}

\subsection{General Framework for Fractionally Dissipative PDEs.}  For $s \in (0,2]$, we consider the system of fractionally dissipative partial differential equations
\be \la{abs}
\pa_t v + c(-\Delta)^{\fr{s}{2}} v = f(v, \tilde{v})
\ee with initial data
\be 
v(x,0) = v_0(x)
\ee on the $d$-dimensional torus $\TT^d$ with periodic boundary conditions, where $v = (v_1, \dots, v_n)$ is an $n$-dimensional vector field depending on a spatial variable $x \in \TT^d$ and a time variable $t \in [0,\infty)$, $c$ is a positive constant, $\tilde{v}$ is a given smooth time-independent function that depends only on $x$, and $f(v, \tilde{v}) = (f_1(v, \tilde{v}), \dots, f_n(v, \tilde{v}))$ is an $n$-dimensional vector field consisting of linear or nonlinear operators $f_1, \dots, f_n$ in $v$. Suppose  \eqref{abs} has a unique global-in-time smooth solution denoted by $v(x,t)$. We establish a bootstrapping framework to address the behavior of higher-order Sobolev norms $v$ in the presence and absence of $\tilde{v}$. We present a proof based on induction, energy estimates, fractional product inequalities and uniform Gronwall lemmas that yields bounds of the form
\be \label{desb}
\|(-\Delta)^{\fr{k}{2}} v(t) \|_{L^2}^2 + \int\limits_{t}^{t+1} \|(-\Delta)^{\fr{k}{2} + \frac{s}{4}} v(s) \|_{L^2}^2 ds \le C_k^1 e^{-c_kt} + C_k^2 
\ee for any $k \ge 1$, where 
\begin{enumerate}
    \item $C_k^1$ is a positive constant depending only on $k$, the size of the initial data $v_0$ in $H^k$, and $\tilde{v}$;
    \item $c_k$ is a positive constant depending only on $k$ and the diffusion coefficient $c$ in \eqref{abs};
    \item $C_k^2$ is a nonnegative constant depending only on $k$ and $\tilde{v}$ such that $C_k^2 = 0 $ when $\tilde{v} = 0$.
    \end{enumerate}
In other words, if $\tilde{v} =0$, then the solution $v$ decays exponentially fast in time to $0$ in all Sobolev norms whereas in the presence of added body functions $\tilde{v}$, the Sobolev norms of $v$ depends only on the size of $\tilde{v}$ for large times. 

In order to obtain such estimates, we apply an induction argument where the base step amounts to having 
\be \la{ind1}
\|\na v(t)\|_{L^2}^2 + \int_{t}^{t+1} \|(-\Delta)^{\fr{s}{4} + \fr{1}{2}} v(s)\|_{L^2}^2 ds \le C_1^2 e^{-c_1 t} + C_1^2
\ee for all $t \ge 0$. Obtaining this first level of regularity could be challenging depending on the structure of the nonlinearity described by $f(v, \tilde{v})$. Assuming that $\eqref{ind1}$ is established and that the inductive hypothesis 
\be \la{ind3}
\|(-\Delta)^{\fr{k-1}{2}} v(t)\|_{L^2}^2  + \int_{t}^{t+1} \|(-\Delta)^{\fr{k-1}{2} + \frac{s}{4}} v(s)\|_{L^2}^2 ds  \le C_{k-1}^1 e^{-c_{k-1}t} + C_{k-1}^2
\ee  holds for all $t \ge 0$,  we prove that the bound \eqref{ind3} holds at the $(k+1)$-th stage. Indeed, we take the $L^2(\TT^d)$ inner product of the equation \eqref{abs} obeyed by $v$ with $(-\Delta)^k v$ and we get the energy equality
\be \la{hhh}
\fr{1}{2} \fr{d}{dt} \|(-\Delta)^{\fr{k}{2}} v\|_{L^2}^2 + c\|(-\Delta)^{\fr{k}{2} + \frac{s}{4}} v\|_{L^2}^2 = \int_{\TT^d} (-\Delta)^{\frac{k}{2} - \frac{s}{4}} f(v, \tilde{v}) \cdot  (-\Delta)^{\frac{k}{2} + \frac{s}{4}} v dx.
\ee We seek estimates of the spacial integral involving the operator $f$ of the form 
\be \la{nonlinearestimate}
\beg{aligned}
&\left|\int_{\TT^d} (-\Delta)^{\frac{k}{2} - \frac{s}{4}}f(v, \tilde{v})  \cdot (-\Delta)^{\frac{k}{2} + \frac{s}{4}} v dx \right|
\le \fr{c}{2}\|(-\Delta)^{\fr{k}{2} + \frac{s}{4}} v\|_{L^2}^2 
\\&\quad+ \left(\Gamma_{1,k} \|(-\Delta)^{\fr{k-1}{2}} v(t)\|_{L^2}^{p_1} + \Gamma_{2,k} (\tilde{v}) \right) \|(-\Delta)^{\fr{k-1}{2} + \frac{s}{4}} v(s)\|_{L^2}^2
+ \left(\Gamma_{3,k} + \Gamma_{4,k} (\tilde{v})\right) \|(-\Delta)^{\fr{k-1}{2}} v(t)\|_{L^2}^{p_2}
\end{aligned}
\ee where 
\begin{enumerate}
\item $\Gamma_{1,k}$ and $\Gamma_{3,k}$ are positive constants depending on $k$;
\item $\Gamma_{2,k}$ and $\Gamma_{4,k}$ are nonnegative constants depending only on $k$ and $\tilde{v}$ such that $\Gamma_{2,k} = \Gamma_{4,k} = 0$ when $\tilde{v} = 0$;
\item $p_1 \ge 0$ and $p_2 \ge 0$ are some nonnegative powers.
    \end{enumerate}
In this case, the energy evolution \eqref{hhh} gives rise to the inequality 
\beg{align} \la{ind5}
&\fr{d}{dt}\|(-\Delta)^{\fr{k}{2}} v\|_{L^2}^2 + c\|(-\Delta)^{\fr{k}{2} + \frac{s}{4}} v\|_{L^2}^2 \nonumber
\\&\le 2 \left(\Gamma_{1,k} \|(-\Delta)^{\fr{k-1}{2}} v\|_{L^2}^{p_1} + \Gamma_{2,k}(\tilde{v})\right)  \|(-\Delta)^{\fr{k-1}{2} + \frac{s}{4}} v\|_{L^2}^2
+ 2 \left(\Gamma_{3,k} + \Gamma_{4,k}(\tilde{v})\right) \|(-\Delta)^{\fr{k-1}{2}} v\|_{L^2}^{p_2}.
\end{align}
In view of the induction hypothesis \eqref{ind3}, the conditions of Proposition \ref{localinteg} are satisfied. Consequently, an application of the uniform Gronwall Lemma \ref{exp} gives the desired bound
\be 
\|(-\Delta)^{\fr{k}{2}} v(t) \|_{L^2}^2 + \int_{t}^{t+1} \|(-\Delta)^{\fr{k}{2} + \frac{s}{4}}v(s) \|_{L^2}^2 ds  \le C_{k}^1e^{-c_k t} + C_k^2.
\ee In summary, if the $H^1$ estimate \eqref{ind1} and the nonlinear estimate \eqref{nonlinearestimate} holds, then one obtain the bounds \eqref{desb} describing the long-time behavior of solutions to \eqref{abs}. We apply this general framework to address the asymptotic behavior of solutions to the NPD system.

\subsection{Higher-order Sobolev Regularity of Solutions to the NPD System.} In this subsection, we study the smoothness and asymptotic behavior of solutions to the NPD system in the presence and absence of body charges $\tilde{\rho}$. We need the following lemma:

\beg{lem} Let $g \in H^3$ be a mean-free function. Suppose $\Psi$ is a solution of the 3D Poisson equation 
\be 
-\Delta \Psi = g
\ee with periodic boundary conditions. Let $v = - \PP(g \na \Psi)$. Then the following interpolation inequalities 
\be \la{poissonineq}
\|\Delta u\|_{L^3}
\le C\|g\|_{H^1}\|g\|_{H^2}^{\frac{1}{2}} \|g\|_{H^3}^{\frac{1}{2}}
\ee and 
\be \la{poissonineq2}
\|\na u\|_{L^2}
\le C\|g\|_{H^1}^2
\ee hold.
\end{lem}

\beg{proof}
As $\Delta$ and $\PP$ commutes and $\PP$ is bounded on $L^3$, we have 
\be 
\|\Delta u\|_{L^3} \le C\|\Delta (g \na \Psi)\|_{L^3}
\le C\|\Delta g \na \Psi\|_{L^3} + C\|\na g \cdot \na \na \Psi\|_{L^3} + C\|g \Delta \na \Psi\|_{L^3}.
\ee Applications of H\"older's inequality yield
\be 
\|\Delta u\|_{L^3}
\le C\|\Delta g\|_{L^3} \|\na \Psi\|_{L^{\infty}} + C\|\na g\|_{L^6} \|\na \na \Psi\|_{L^6} + C\|g\|_{L^6} \|\Delta \na \Psi\|_{L^6}.
\ee By making use of standard 3D continuous Sobolev embeddings, we have
\be 
\|\na \Psi\|_{L^{\infty}} 
\le C\|\na \Psi\|_{H^2}
\le C\|\Psi\|_{H^3}
\le C\|(-\Delta)^{-1} g\|_{H^3}
\le C\|g\|_{H^1},
\ee 
\be 
\|\na \na \Psi\|_{L^6}
\le C\|\na \na \Psi\|_{H^1}
\le C\|\Psi\|_{H^3}
= C\|(-\Delta)^{-1} g\|_{H^3}
\le C\|g\|_{H^1},
\ee and 
\be 
\|\Delta \na \Psi\|_{L^6}
\le C\|\Delta \na \Psi\|_{H^1}
\le C\|\Psi\|_{H^4}
\le C\|(-\Delta)^{-1} g\|_{H^4}
\le C\|g\|_{H^2}.
\ee 
Consequently, it follows that 
\be 
\beg{aligned}
\|\Delta u\|_{L^3}
&\le C\|\Delta g\|_{L^2}^{\frac{1}{2}}\|\na \Delta g\|_{L^2}^{\frac{1}{2}} \|g\|_{H^1}
+ C\|\Delta g\|_{L^2} \|g\|_{H^1}
+ C\|\na g\|_{L^2} \|g\|_{H^2}
\\&\le C\|g\|_{H^2}^{\frac{1}{2}}\|g\|_{H^3}^{\frac{1}{2}} \|g\|_{H^1}
+ C\|g\|_{H^2}\|g\|_{H^1}.
\end{aligned} 
\ee Finally, we apply the inequality $\|g\|_{H^2}^{\frac{1}{2}} \le \|g\|_{H^2}^{\frac{1}{2}}\|g\|_{H^3}^{\frac{1}{2}}$, and deduce \eqref{poissonineq}. As for \eqref{poissonineq2}, we have 
\be 
\beg{aligned}
\|\na u\|_{L^2}
\le C\|\na (g \na \psi)\|_{L^2}
\le C\|\na g\|_{L^2} \|\na \Psi\|_{L^{\infty}}
+ C\|g\|_{L^3} \|\na \na \Psi\|_{L^6}
\le C\|g\|_{H^1}^2.
\end{aligned}
\ee 
\end{proof}

\begin{prop} \la{mpfora} For any time $t \ge 0$, it holds that 
\be \label{desb}
\sum\limits_{i=1}^{N} \|(-\Delta)^{\fr{k}{2}} c_i(t) \|_{L^2}^2 +\int_{t}^{t+1} \sum\limits_{i=1}^{N} \|(-\Delta)^{\fr{k+1}{2}} c_i(s) \|_{L^2}^2 ds \le C_k^1 e^{-c_kt} + C_k^2 
\ee for any integer $k \ge 1$, where $C_k^1$ is a positive constant depending only on $k$, the size of the initial data $v_0$ in $H^k$, and $\tilde{\rho}$, $c_k$ is a positive constant depending only on $k$ and the diffusivity $D$, and  $C_k^2$ is a nonnegative constant depending only on $k$ and $\tilde{\rho}$ such that $C_k^2 = 0 $ when $\tilde{\rho} = 0$.
\end{prop}

\beg{proof} The proof is based on the general framework introduced in the previous subsection where $v = (c_1, \dots, c_N)$, $\tilde{v} = \tilde{\rho}$, $s =2$, $c=D$, and $f(v, \tilde{v}) = (- u \cdot \na c_1 +Dz_1 \na \cdot (c_1 \na \Phi), \dots,- u \cdot \na c_N +Dz_N \na \cdot (c_N \na \Phi) )$. When $k=1$, the result follows from Theorem \ref{T1} and yields the $H^1$ estimate corresponding to \eqref{ind1}. We prove that the nonlinear estimate \eqref{nonlinearestimate} holds. We fix an integer $k \ge 3$. Using the fact that $H^{k-1}$ is a Banach Algebra, the boundedness of the Leray projector on $L^2$, and elliptic estimates, we have  
\be 
\beg{aligned}
&\left|\sum\limits_{i=1}^{N} \int_{\TT^3} (-\Delta)^{\frac{k-1}{2}} (u \cdot \na c_i) (-\Delta)^{\frac{k+1}{2}} c_i dx\right|
\le \frac{D}{4} \sum\limits_{i=1}^{N}\|(-\Delta)^{\frac{k+1}{2}} c_i\|_{L^2}^2
+ C\sum\limits_{i=1}^{N}\|(-\Delta)^{\frac{k-1}{2}} (u \cdot \na c_i)\|_{L^2}^2
\\&\quad\quad\le \frac{D}{4} \sum\limits_{i=1}^{N}\|(-\Delta)^{\frac{k+1}{2}} c_i\|_{L^2}^2 + C \sum\limits_{i=1}^{N} \|(-\Delta)^{\frac{k-1}{2}}u\|_{L^2}^2 \|(-\Delta)^{\frac{k}{2}} c_i\|_{L^2}^2
\\&\quad\quad\le \frac{D}{4} \sum\limits_{i=1}^{N}\|(-\Delta)^{\frac{k+1}{2}} c_i\|_{L^2}^2 + C\sum\limits_{i=1}^{N}\|(-\Delta)^{\frac{k-1}{2}}\PP ((\rho +\tilde{\rho}) \na \Phi)\|_{L^2}^2 \|(-\Delta)^{\frac{k}{2}} c_i\|_{L^2}^2
\\&\quad\quad\le  \frac{D}{4} \sum\limits_{i=1}^{N} \|(-\Delta)^{\frac{k+1}{2}} c_i\|_{L^2}^2 + C\|(-\Delta)^{\frac{k-1}{2}} (\rho +\tilde{\rho})\|_{L^2}^2 (-\Delta)^{\frac{k-1}{2}} \na \Phi\|_{L^2}^2  \sum\limits_{i=1}^{N} \|(-\Delta)^{\frac{k}{2}} c_i\|_{L^2}^2
\\&\quad\quad\le \frac{D}{4} \sum\limits_{i=1}^{N} \|(-\Delta)^{\frac{k+1}{2}} c_i\|_{L^2}^2 + C\left(\sum\limits_{i=1}^{N} \left(\|(-\Delta)^{\frac{k-1}{2}} c_i\|_{L^2}^2\right)^2 +  C\|\tilde{\rho}\|_{H^{k-1}}^4\right) \sum\limits_{i=1}^{N} \|(-\Delta)^{\frac{k}{2}} c_i\|_{L^2}^2.
\end{aligned}
\ee If $k=2$, we make use of the cancellation law 
\be 
\int_{\TT^3} (u \cdot \na \Delta c_i) \Delta c_i dx = 0
\ee together with the interpolation inequalities \eqref{poissonineq} and \eqref{poissonineq2} to estimate
\be 
\beg{aligned}
&\left|\sum\limits_{i=1}^{N} \int_{\TT^3} (-\Delta)^{\frac{1}{2}} (u \cdot \na c_i) (-\Delta)^{\frac{3}{2}} c_i dx \right|
= \left|\sum\limits_{i=1}^{N} \int_{\TT^3} (-\Delta)(u \cdot \na c_i) (-\Delta) c_i dx \right|
\\&\quad\quad\le C\sum\limits_{i=1}^{N} \left[\|\Delta u\|_{L^3} \|\na c_i\|_{L^2} \|\Delta c_i\|_{L^6} + \|\na u\|_{L^2} \|\na \na c_i\|_{L^3} \|\Delta c_i\|_{L^6} \right]
\\&\quad\quad\le C\sum\limits_{i=1}^{N} \left(\|\rho + \tilde{\rho}\|_{H^1} \|\rho + \tilde{\rho}\|_{H^2}^{\frac{1}{2}} \|\rho + \tilde{\rho}\|_{H^3}^{\frac{1}{2}} \right)\|\na c_i\|_{L^2}  \|(-\Delta)^{\frac{3}{2}}c_i\|_{L^2}  \\&\quad\quad\quad\quad+C\sum\limits_{i=1}^{N}   \|\rho + \tilde{\rho}\|_{H^1}^2 \|\Delta c_i\|_{L^2}^{\frac{1}{2}} \|(-\Delta)^{\frac{3}{2}}c_i \|_{L^2}^{\frac{3}{2}}
\\&\le \frac{D}{4}\sum\limits_{i=1}^{N}\|(-\Delta)^{\frac{3}{2}}c_i\|_{L^2}^2 + C\left( \left(\sum\limits_{i=1}^{N} \|\na c_i \|_{L^2}^2 \right)^4 +\|\tilde{\rho}\|_{H^3}^8 + 1 \right)\sum\limits_{i=1}^{N} \|\Delta c_i\|_{L^2}^2.
\end{aligned}
\ee 
As for the electromigration term, since $H^{k}$ is a Banach Algebra for $k \ge 2$, we have 
\be 
\beg{aligned}
&\left|\sum\limits_{i=1}^{N} \int_{\TT^3} (-\Delta)^{\frac{k-1}{2}} \na \cdot (c_i \na \Phi) (-\Delta)^{\frac{k+1}{2}} c_i dx\right|
\le \frac{D}{4} \sum\limits_{i=1}^{N} \|(-\Delta)^{\frac{k+1}{2}} c_i\|_{L^2}^2
+ C \sum\limits_{i=1}^{N} \|(-\Delta)^{\frac{k}{2}} (c_i \na \Phi)\|_{L^2}^2
\\&\quad\quad\le \frac{D}{4} \sum\limits_{i=1}^{N} \|(-\Delta)^{\frac{k+1}{2}} c_i\|_{L^2}^2
+ C \|(-\Delta)^{\frac{k}{2}} \na \Phi \|_{L^2}^2 \sum\limits_{i=1}^{N} \|(-\Delta)^{\frac{k}{2}} c_i \|_{L^2}^2
\\&\quad\quad\le \frac{D}{4} \sum\limits_{i=1}^{N} \|(-\Delta)^{\frac{k+1}{2}} c_i\|_{L^2}^2 + C\left(\sum\limits_{i=1}^{N} \|(-\Delta)^{\frac{k-1}{2}} c_i\|_{L^2}^2 +  C\|\tilde{\rho}\|_{H^{k-1}}^2\right) \sum\limits_{i=1}^{N} \|(-\Delta)^{\frac{k}{2}} c_i\|_{L^2}^2.
\end{aligned}
\ee Therefore, the nonlinear estimate \eqref{nonlinearestimate} is fulfilled for the 3D NPD system. 
\end{proof}

\section{Properties of the Solution Map} \la{s5}

Let $\mathcal{V}$  be the subset of $H^1 \times \dots \times  H^1$ consisting of vectors $(\xi_1, \dots, \xi_N)$ such that:
\begin{enumerate}
\item $\xi_1, \dots, \xi_N$ are nonnegative functions almost everywhere;
\item The density $R = \sum\limits_{i=1}^{N} z_i \xi_i$ obeys $\int_{\TT^3} (R + \tilde{\rho}) dx = 0$.
\end{enumerate} We define the solution map 
\begin{equation}
\mathcal{S}(t): \mathcal{V} \mapsto \mathcal{V}
\end{equation}  
corresponding to the 3D periodic NPD system  by 
\begin{equation}
\mathcal{S}(t) (c_1(0), \dots, c_N(0)) = (c_1(t), \dots, c_N(t)).
\end{equation} 
For each $M > 0$, we consider the subset $\mathcal{V}_M$ of $\mathcal{V}$ consisting of vectors $v = (\xi_1, \dots, \xi_N) \in \mathcal{V}$ obeying $\int_{\TT^3} v(x) dx \le M$.  

\beg{prop} \la{exabs}
Fix $M > 0$ and an integer $k \ge 2$. There exists a radius $R > 0$ depending only on the density $\tilde{\rho}$, $M$ and $k$ such that for any 
$(c_1(0), \dots, c_N(0)) \in \mathcal{V}_M$, there exists $t_0 > 0$ depending on 
the $H^1$ norm of $(c_1(0), \dots, c_N(0))$ such that 
\begin{equation}
\mathcal{S}(t)(c_1(0), \dots, c_N(0)) \in \mathcal{B}_{R,M,k} 
= \left\{(c_1, \dots, c_N) \in \mathcal{V}_M \cap (H^k \times \dots \times H^k): \sum\limits_{i=1}^{N} \|c_i - \bar{c}_i\|_{H^k} \leq R \right\}
\end{equation} 
for all $t \geq t_0$.
\end{prop}

\begin{proof}
The proof follows from Proposition \ref{mpfora}.
\end{proof}

\beg{rem} \la{propabs} The ball $\mathcal{B}_{R, M, k}$ is compact in $\mathcal{V}$ for any $k \ge 2$, a fact that follows from the uniform boundedness of the spatial averages of the ionic concentrations by $M$, the Rellich-Kondrachov theorem, and the Banach-Alaoglu theorem.   Also, $\mathcal{B}_{R, M, k}$ is convex as $\mathcal{V}_M$ is convex, and so it is connected. 
\end{rem}

The solution map is instantaneously Lipschitz continuous on $\mathcal{V}$ as stated in the following proposition:

\beg{prop} \la{contsol}
Let $w_1^0 = (c_1^1(0), \dots, c_N^1(0)), w_2^0 = (c_2^1(0), \dots, c_N^2(0)) \in \mathcal{V}$. It holds that 
\be \label{ContV}
\|\mathcal{S}(t)w_1^0 - \mathcal{S}(t)w_2^0\|_{\mathcal{V}}^2 \leq CK(t) \|w_1^0 - w_2^0\|_{\mathcal{V}}^2,
\end{equation} for any $t \ge 0$, where 
\be 
K(t) = e^{\int_{0}^{t}   \left(1 + \|\tilde{\rho}\|_{H^1}^2 +   \| \mathcal{S}(t)w_1^0 \|_{H^1}^{2} + \|\mathcal{S}(t)w_2^0\|_{H^1}^{2}  \right)\left( \| \mathcal{S}(t)w_1^0 \|_{H^2}^{2} + \|\mathcal{S}(t)w_2^0\|_{H^2}^{2} + \|\tilde{\rho}\|_{H^2}^{2} + 1\right)  ds }.
\ee 
\end{prop}

The proof of Proposition \ref{contsol} is given in Appendix \ref{AP2}. 

\section{Backward Uniqueness}
\la{s6}

In this section, we generalize an argument of \cite{constantin2015long} and present a general framework to study the backward uniqueness of solutions to fractionally dissipative PDEs, and we apply it to the NPD system. 

\subsection{General Framework for Backward Uniqueness.}  For $s \in (0,2]$, we consider the system of fractionally dissipative partial differential equation
\be \la{absu}
\pa_t v + c(-\Delta)^{\fr{s}{2}} v = f(v, \tilde{v})
\ee with initial data
\be 
v(x,0) = v_0(x)
\ee on the $d$-dimensional torus $\TT^d$ with periodic boundary conditions, where $v = (v_1, \dots, v_n)$ is an $n$-dimensional vector field depending on a spatial variable $x \in \TT^d$ and a time variable $t \in [0,\infty)$, $c$ is a positive constant, $\tilde{v}$ is a given smooth time-independent function that depends only on $x$, and $f(v, \tilde{v}) = (f_1(v, \tilde{v}), \dots, f_n(v, \tilde{v}))$ is an $n$-dimensional vector field consisting of linear or nonlinear operators $f_1, \dots, f_n$ in $v$. Suppose  \eqref{absu} has a unique global-in-time smooth solution denoted by $v(x,t)$. 

\beg{Thm} \la{frameback} Let $v_1$ and $v_2$ be two smooth solutions to \eqref{absu} corresponding to two initial data $v_1(0)$ and $v_2(0)$ respectively. Suppose there is a time $T$ such that $v_1(T) = v_2(T)$. Denote by $v$ the difference $v_1 - v_2$. If 
\be \la{baun4}
\|f(v_1, \tilde{v}) - f(v_2, \tilde{v})\|_{L^2(\TT^d)}^2 \le K_1(t) \|v\|_{L^2(\TT^d)}^2 + K_2(t) \|(-\Delta)^{\frac{s}{4}} v\|_{L^2(\TT^d)}^2,
\ee on $[0,T]$ with 
\be 
\int_{0}^{T} (K_1(t) + K_2(t)) dt < \infty,
\ee then $v_1(0) = v_2(0)$. 
\end{Thm}

\beg{proof}
Without loss of generality, suppose $v \ne 0$ on $[0,T)$. Let $E_0 = \|v\|_{L^2(\TT^d)}^2$, $E_1 = \|(-\Delta)^{\frac{s}{4}} v\|_{L^2(\TT^d)}^2$ and $Y(t) = - \log E_0$. We note that $Y(t)$ is bounded from below and blows up at $t$ approaches $T^-$. We show that 
\be \la{baun1}
\frac{d}{dt}\frac{E_1}{E_0} \le C_1(t) \frac{E_1}{E_0} + C_2(t)
\ee and 
\be \la{baun2}
\frac{d}{dt} Y \le C_3(t)\frac{E_1}{E_0} + C_4(t)
\ee for some functions $C_1, \dots, C_4$ obeying 
\be 
\int_{0}^{T} \left(C_1 + C_2 + C_3 + C_4\right) dt <  \infty.
\ee This implies that $Y \in L^{\infty}(0,T)$, contradicting the blow up near $T$. The energies $E_0$ and $E_1$ evolve according to 
\be 
\frac{1}{2}\frac{d}{dt} E_0 + cE_1 = (f(v_1, \tilde{v}) - f(v_2, \tilde{v}), v)_{L^2(\TT^d)}
\ee and 
\be 
\frac{1}{2}\frac{d}{dt}  E_1 + c\|(-\Delta)^{\frac{s}{2}} v\|_{L^2(\TT^d)}^2 = (f(v_1, \tilde{v}) - f(v_2, \tilde{v}), (-\Delta)^{\frac{s}{2}} v)_{L^2(\TT^d)}.
\ee Since $Y$ obeys
\be 
\frac{1}{2} \frac{dY}{dt} = -\frac{1}{2} \frac{d}{dt} \log E_0 = c\frac{E_1}{E_0} - \frac{(f(v_1, \tilde{v}) - f(v_2, \tilde{v}), v)_{L^2(\TT^d)}}{E_0}
\ee and 
\be 
\beg{aligned}
|(f(v_1, \tilde{v}) - f(v_2, \tilde{v}), v)_{L^2(\TT^d)} |
&\le \frac{1}{2} \|f(v_1, \tilde{v}) - f(v_2, \tilde{v}) \|_{L^2(\TT^d)}^2 + \frac{1}{2} \|v\|_{L^2(\TT^d)}^2
\\&\le \frac{1}{2} (K_1(t) E_0 + K_2(t)E_1 + E_0),
\end{aligned}
\ee we infer that the differential inequality \eqref{baun2} holds. Now we turn our attention to the time evolution of $\frac{E_1}{E_0}$ and derive \eqref{baun1}. We have 
\be 
\beg{aligned}
\frac{1}{2} \frac{d}{dt} \frac{E_1}{E_0} 
&= \frac{1}{2E_0} \frac{dE_1}{dt}  + \frac{1}{2}E_1 \frac{d}{dt} \frac{1}{E_0}
\\&\le -\frac{c}{E_0} \|(-\Delta)^{\frac{s}{2}} v\|_{L^2}^2 + \frac{cE_1^2}{E_0^2} + \frac{1}{E_0} (f(v_1, \tilde{v}) - f(v_2, \tilde{v}), ((-\Delta)^{\frac{s}{2}} - E_1E_0^{-1})v)_{L^2(\TT^d)}.
\end{aligned}
\ee We compute the following inner product
\be 
\beg{aligned}
 &(((-\Delta)^{\frac{s}{2}} - E_1E_0^{-1})v, ((-\Delta)^{\frac{s}{2}} - E_1E_0^{-1})v)_{L^2(\TT^d)}
\\&\quad\quad= \|(-\Delta)^{\frac{s}{2}} v\|_{L^2(\TT^d)}^2 + E_1^2E_0^{-2} \|v\|_{L^2(\TT^d)}^2 -2E_1E_0^{-1} ((-\Delta)^{\frac{s}{2}}v, v)_{L^2(\TT^d)}
\\&\quad\quad= \|(-\Delta)^{\frac{s}{2}} v\|_{L^2(\TT^d)}^2 + E_1^2E_0^{-1} -2E_1^2E_0^{-1}
\\&\quad\quad= \|(-\Delta)^{\frac{s}{2}} v\|_{L^2(\TT^d)}^2 - E_1^2E_0^{-1},
\end{aligned}
\ee and deduce the identity
\be 
 -\frac{c}{E_0} \|(-\Delta)^{\frac{s}{2}} v\|_{L^2}^2 + \frac{cE_1^2}{E_0^2} = -\frac{c}{E_0} \|((-\Delta)^{\frac{s}{2}} - E_1E_0^{-1}) v\|_{L^2(\TT^d)}^2.
\ee Consequently, we obtain 
\be 
\frac{1}{2} \frac{d}{dt} \frac{E_1}{E_0} + \frac{1}{2}\frac{c}{E_0} \|((-\Delta)^{\frac{s}{2}} - E_1E_0^{-1}) v\|_{L^2(\TT^d)}^2 \le CE_0^{-1} \| f(v_1, \tilde{v}) - f(v_2, \tilde{v})\|_{L^2(\TT^d)}^2,
\ee from which we deduce \eqref{baun1} after making use of \eqref{baun4}.
\end{proof}

\subsection{Injectivity of the Solution Map} In this subsection, we use the aforementioned general framework to prove the injectivity of the solution map.

\beg{prop} \la{injsol}
Let $w_1 = (c_1^1, \dots, c_N^1)$ and $w_2 = (c_1^2, \dots, c_N^2)$ be solutions to the NPD system corresponding to two smooth initial data $w_1(0) = (c_1^1(0), \dots, c_N^1(0))  \in \mathcal{V}$ and $w_2(0) = (c_1^2(0), \dots, c_N^2(0)) \in \mathcal{V}$ respectively. Suppose there is a time $T$ such that $w_1(T) = w_2(T)$. Then $w_1(0) = w_2(0)$. 
\end{prop}

\beg{proof}
We define the differences
\be 
w= w_1 - w_2, c_i = c_i^1 - c_i^2, u = u_1 - u_2, \rho = \rho_1 - \rho_2, \Phi = \Phi_1 - \Phi_2
\ee where 
\be 
\beg{aligned}
&\rho_1 = \sum\limits_{i=1}^{N} z_i c_i^1, \rho_2 = \sum\limits_{i=1}^{N} z_i c_i^2, -\Delta \Phi_1 = \rho_1 + \tilde{\rho}, -\Delta \Phi_2 = \rho_2 + \tilde{\rho}, 
\\&u_1 = -\PP((\rho_1 + \tilde{\rho}) \na \Phi_1), u_2 = -\PP((\rho_2 + \tilde{\rho}) \na \Phi_2).
\end{aligned}
\ee 

{\bf{Step 1. Velocity Estimates.}} The difference $u$ obeys the constitutive relation 
\be 
u = -\PP(\rho \na \Phi_1) - \PP((\rho_2 + \tilde{\rho})\na \Phi)
\ee where $\Phi$ solves the Poisson equation 
\be 
-\Delta \Phi =\rho.
\ee Due to the boundedness of the Leray projector on $L^3$, we have
\be
\|u\|_{L^3}^2 
\le C\|\rho \na \Phi_1\|_{L^3}^2 + C\|(\rho_2 + \tilde{\rho}) \na \Phi)\|_{L^3}^2
\le C\|\rho\|_{L^6}^2 \|\na \Phi_1\|_{L^6}^2 + C\|\rho_2 + \tilde{\rho}\|_{L^6}^2 \|\na \Phi\|_{L^6}^2,
\ee which boils down to 
\be 
\|u\|_{L^3}^2
\le C\|\na \rho\|_{L^2}^2 \|\rho_1 + \tilde{\rho}\|_{L^2}^2 + C\|\rho_2 + \tilde{\rho}\|_{H^1}^2 \|\rho\|_{L^2}^2
\ee after making use of elliptic estimates and the 3D continuous Sobolev embedding of $H^1$ in $L^6$. Consequently, we infer that 
\be \la{veloest1}
\|u\|_{L^3}^2 \le C \left(\sum\limits_{i=1}^{N} \|c_i^1\|_{H^1}^2 + \|c_i^2\|_{H^1}^2 + \|\tilde{\rho}\|_{H^1}^2 \right) \left(\sum\limits_{i=1}^{N} \|c_i\|_{L^2}^2 + \|\na c_i\|_{L^2}^2 \right).
\ee Since $\PP$ is bounded on $H^2$ and $H^2$ is a Banach Algebra, we have
\be 
\|\Delta u_2\|_{L^2}^2
\le C\|(\rho_2 + \tilde{\rho}) \na \Phi_2\|_{H^2}^2
\le C\|\rho_2 + \tilde{\rho}\|_{H^2}^2 \|\na \Phi_2\|_{H^2}^2
\le C\|\rho_2 + \tilde{\rho}\|_{H^2}^2 \|\rho_2 + \tilde{\rho}\|_{H^1}^2,
\ee giving rise to the estimate 
\be \la{veloest2}
\|\Delta u_2\|_{L^2}^2
\le C\left(\|\tilde{\rho}\|_{H^2}^2 + \sum\limits_{i=1}^{N} \|c_i^2\|_{H^2}^2 \right)\left(\sum\limits_{i=1}^{N} \|c_i^2\|_{H^1}^2 +\|\tilde{\rho}\|_{H^1}^2\right).
\ee 

{\bf{Step 2. Backward Uniqueness.} }For a vector $\omega = (\xi_1, \dots, \xi_N)$, we let $R = \sum\limits_{i=1}^{N} z_i \xi_i$, $\Psi = (-\Delta)^{-1} (R+\tilde{\rho})$, and $u = -\PP((R+\tilde{\rho})\na \Psi)$, and we consider the operator 
\be 
f(\omega, \tilde{\rho}) = (- u \cdot \na \xi_1 + Dz_i \na \cdot (\xi_1 \na \Psi), \dots, - u \cdot \na \xi_N + Dz_N \na \cdot (\xi_N \na \Psi)).
\ee We show that 
\be \la{backuniqueness}
\|f(w_1, \tilde{\rho}) - f(w_2, \tilde{\rho})\|_{L^2}^2 
\le C_1(t) \|w\|_{L^2}^2 +C_2(t) \|\na w\|_{L^2}^2
\ee for some time integrable functions $C_1$ and $C_2$ over $[0,T]$, and we obtain the desired result as a consequence of Theorem \ref{frameback}. Indeed, it follows from \eqref{veloest1} and \eqref{veloest2} that 
\be 
\beg{aligned}
&\sum\limits_{i=1}^{N} \|u_1 \cdot \na c_i^1 - u_2 \cdot \na c_i^2\|_{L^2}^2
= \sum\limits_{i=1}^{N} \|u \cdot \na c_i^1 + u_2 \cdot \na c_i\|_{L^2}^2
\\&\quad\quad\le \sum\limits_{i=1}^{N} \|u\|_{L^3}^2 \|\na c_i^1\|_{L^6}^2 + \|u_2\|_{L^{\infty}}^2 \|\na c_i\|_{L^2}^2
\le C\|u\|_{L^3}^2 \sum\limits_{i=1}^{N} \|\Delta c_i^1\|_{L^2}^2 + C\|\Delta u_2\|_{L^2}^2 \sum\limits_{i=1}^{N} \|\na c_i\|_{L^2}^2
\\&\quad\quad\le C\left(\|\tilde{\rho}\|_{H^2}^2 + \sum\limits_{i=1}^{N} \|c_i^1\|_{H^1}^2 + \|c_i^2\|_{H^1}^2 \right) \left(\|\tilde{\rho}\|_{H^2}^2 +\sum\limits_{i=1}^{N} \|c_i^1\|_{H^2}^2 + \|c_i^2\|_{H^2}^2 \right) \sum\limits_{i=1}^{N} \left(\|c_i\|_{L^2}^2 + \|\na c_i\|_{L^2}^2\right).
\end{aligned}
\ee As for the electromigration difference, we have
\be 
\beg{aligned}
&\sum\limits_{i=1}^{N} \|\na \cdot (c_i^1 \na \Phi_1) - \na \cdot (c_i^2 \na \Phi_2)\|_{L^2}^2
= \sum\limits_{i=1}^{N} \|\na \cdot (c_i \na \Phi_1) + \na \cdot (c_i^2 \na \Phi)\|_{L^2}^2
\\&\quad= \sum\limits_{i=1}^{N} \|\na c_i \cdot \na \Phi_1 + c_i \Delta \Phi_1 + \na c_i^2 \cdot \na \Phi + c_i^2 \Delta \Phi\|_{L^2}^2
\\&\quad\le \sum\limits_{i=1}^{N} \|\na c_i\|_{L^2}^2 \|\na \Phi_1\|_{L^{\infty}}^2 + \|c_i\|_{L^2}^2 \|\Delta \Phi_1\|_{L^{\infty}}^2 + \|\na c_i^2\|_{L^2}^2 \|\na \Phi\|_{L^{\infty}}^2 + \|c_i^2\|_{L^{\infty}}^2\|\Delta \Phi\|_{L^2}^2
\\&\quad\le C\|\na (\rho_1+\tilde{\rho})\|_{L^2}^2 \sum\limits_{i=1}^{N} \|\na c_i\|_{L^2}^2 + \|\Delta (\rho_1 + \tilde{\rho})\|_{L^2}^2 \sum\limits_{i=1}^{N}\|c_i\|_{L^2}^2 + \|\na \rho\|_{L^2}^2 \sum\limits_{i=1}^{N}\|\na c_i^2\|_{L^2}^2 + C\|\rho\|_{L^2}^2 \sum\limits_{i=1}^{N} \|c_i^2\|_{H^2}^2  
\\&\quad\le C\left(\|\tilde{\rho}\|_{H^2}^2 + \sum\limits_{i=1}^{N} (\|c_i^1\|_{H^2}^2 + \|c_i^2\|_{H^2}^2) \right) \sum\limits_{i=1}^{N} \left(\|c_i\|_{L^2}^2 + \|\na c_i\|_{L^2}^2\right).
\end{aligned}
\ee In the above estimates, we appealed to elliptic estimates and continuous Sobolev embeddings. Since the concentrations belong to the Lebesgue spaces $L^{\infty}(0,T; H^1)$ and $L^2(0,T; H^2)$ we obtain the desired bound \eqref{backuniqueness}. 
\end{proof}

\section{Global Attractor} \la{s7}

In this section, we address the existence of a finite-dimensional global attractor.

\begin{thm}\la{attp}  There exist a time $\tilde T > 0$ and a radius ${R} > 0$ depending only on the body charges $\tilde{\rho}$ and $M$  such that the ball
\begin{equation}
\mathcal{B}_{M} 
= \left\{(c_1, \dots, c_N) \in \mathcal{V}_M \cap (H^2 \times \dots \times H^2): \sum\limits_{i=1}^{N} \|c_i - \bar{c}_i\|_{H^2} \leq R \right\}
\end{equation}  obeys $\mathcal{S}(t) {\mathcal{B}}_M \subset  {\mathcal{B}}_M$ for all $t \ge \tilde T$. 
Moreover, the set \begin{equation} 
{X}_M = \bigcap_{t > 0} S(t) {\mathcal{B}}_M.
\end{equation}
satisfies the following properties: 
\begin{enumerate}
\item[(a)] $X_M$ is compact in $\mathcal V$.
\item[(b)] $S(t) X = X$ for all $t \geq 0$.
\item[(c)] If $Z$ is bounded in $\mathcal V_M$  in the norm of of $\mathcal V$, and $S(t)Z = Z$ for all $t \geq 0$, then $Z \subset  X_M$. 
\item[(d)] For every $w_0 \in \mathcal V_M,$
$\lim\limits_{t \to \infty} dist_{\mathcal V} (S(t)w_0,  X_M) = 0$.
\item[(e)] $X_M$ is connected.
\item[(f)] ${X}_M$ is smooth, that is for every integer $k>0$, there exists a radius $\rho_k$ depending on the body charges $\tilde{\rho}$ and $M$ a ball $B_{\rho_k} \subset H^k \times \dots \times H^k$ such that the attractor ${X}_M \subset B_{\rho_k}$.  
\end{enumerate}
\end{thm}

\beg{proof}
The existence of the absorbing ball $\mathcal{B}_M$ follows from Proposition \ref{exabs}. By making use of the continuity and injectivity of the solution map established in propositions \ref{contsol} and \ref{injsol}  respectively, in addition to the compactness and connectedness of the absorbing ball justified in Remark \ref{propabs}, we can follow the proof of the similar result in \cite{constantin1988navier} to deduce the properties (a)-(e). Finally, the smoothness property (f) of the attractor is a direct consequence of Proposition \ref{exabs}. 
\end{proof}

Fix an integer $n \ge 1$ and let $\mathcal{F}: \Omega \subset \RR^n \rightarrow \mathcal{V}$ be a smooth function. Let $\Sigma_t = \mathcal{S}(t) \mathcal{F}(\Omega)$ and define   the volume element in $\Sigma_t$ to be  
\be 
\left|\frac{\partial}{\partial \alpha_1} \mathcal{S}(t) \mathcal{F}(\alpha) \wedge \dots \wedge \frac{\partial}{\partial \alpha_n} \mathcal{S}(t) \mathcal{F}(\alpha) \right| d\alpha_1 \dots d\alpha_N,
\ee where $d\alpha_1 \dots d\alpha_n$ is the volume element in $\RR^n$ and  $\alpha = (\alpha_1, \dots, \alpha_n) \in \RR^n$. For $j=1,\dots,n $ let 
\be 
\xi^j = (\xi_1^j, \cdots, \xi_N^j):= \frac{\partial}{\alpha_j} \mathcal{S}(t) \mathcal{F}(\alpha), R_j = \sum\limits_{i=1}^{N} z_i \xi_i^j, \Psi_j = (-\Delta)^{-1} R_j.
\ee 
For each $i \in \left\{1, \dots, N\right\}$ and $j \in \left\{1, \dots, n \right\}$, the scalar function $\xi_i^j$ solves the equation 
\be \la{ddde}
\pa_t \xi_i^j - D\Delta \xi_i^j =  \PP(R_j \na \Phi + (\rho + \tilde{\rho})\na\Psi_j) \cdot \na c_i - u \cdot \na \xi_i^j + Dz_i \na \cdot (\xi_{i}^j \na \Phi + c_i \na \Psi_j) 
\ee along $w =  (c_1, \dots , c_N) = \mathcal{S}(t) \mathcal{F}(\alpha)$ with $\rho = \sum\limits_{i=1}^{N} z_i c_i, \Phi= (-\Delta)^{-1}(\rho + \tilde{\rho})$, and $u = - \PP((\rho +\tilde{\rho})\na \Phi)$. For $w_0 = (c_1(0), \dots, c_N(0)) \in X$ we define the  volume
\be 
V_n(t) = \left\|\xi^1 \wedge \dots \wedge \xi^n \right\|_{\bigwedge^n \mathcal{V}}
\ee where $\xi^1, \dots, \xi^n$ are the solutions to \eqref{ddde} along $w(t) = \mathcal{S}(t) w_0$ and $\bigwedge^n \mathcal{V}$ is the $n$-th exterior product of $\mathcal{V}$.  

\beg{prop} \la{dve}
There is a positive time $t_0$ depending only on $R$, an integer $N_0$ depending only on $R$, and a positive constant $c$  such that 
\be 
V_n(t) \le V_n(0) e^{-cn^{2}D t}
\ee for any $t \ge t_0$.
\end{prop}

The proof is based on \cite{constantin1988navier} and is presented in Appendix \ref{AP3} for completion. 

Following the proof of the similar result in \cite{constantin1988navier}, we obtain the following theorem:

\beg{thm}
The attractor $X_M$ has a finite fractal dimension in $\mathcal{V}$, that is there exists a finite real number $\tilde M > 0$ depending on the body charges $\tilde{\rho}$ such that 
\[ 
\lim\sup_{r\to 0}\fr{\log{N_{\mathcal{V}}(r)}}{\log\left(\fr{1}{r}\right)} \le \tilde M 
\] where $N_{\mathcal{V}}(r)$ is the minimal number of balls in $\mathcal{V}$ of radii $r$ needed to cover $X_M$.
\end{thm}

\appendix

\section{Local Well-Posedness} \la{AP1}

In this appendix, we present the proof of Proposition \ref{locso}.

\begin{proof}
The proof is standard and based on Galerkin approximations, energy estimates, and passage to the limit using the Aubin-Lions lemma. We only provide a priori bounds for simplicity. 
The time evolution of $\|\na c_i\|_{L^2}$ is described by 
\be \la{localco}
\frac{1}{2} \frac{d}{dt} \|\na c_i\|_{L^2}^2 
+ D\|\Delta c_i\|_{L^2}^2 
= \int_{\TT^3} u \cdot \na c_i \Delta c_i dx + Dz_i \int_{\TT^3} \na \cdot (c_i \na \Phi) \Delta c_i dx.
\ee We estimate 
\be 
\beg{aligned}
\left|\int_{\TT^3} u \cdot \na c_i \Delta c_i dx \right|
&\le C\|u\|_{L^6}\|\na c_i\|_{L^3}\|\Delta c_i\|_{L^2}
\le C\|(\rho + \tilde{\rho}) \na \Phi\|_{L^6} \|\na c_i\|_{L^2}^{\frac{1}{2}}\|\Delta c_i\|_{L^2}^{\frac{3}{2}} 
\\&\le C\|\rho + \tilde{\rho}\|_{L^6}^4 \|\na \Phi\|_{L^{\infty}}^4 \|\na c_i\|_{L^2}^2 + \frac{D}{4}\|\Delta c_i\|_{L^2}^2
\\&\le C\|\na (\rho + \tilde{\rho})\|_{L^2}^4 \|\na (\rho + \tilde{\rho})\|_{L^2}^4 \|\na c_i\|_{L^2}^2 + \frac{D}{4} \|\Delta c_i\|_{L^2}^2
\\&\le C\left(\sum\limits_{j=1}^{N}\|\na c_j\|_{L^2}^2  + \|\na \tilde{\rho}\|_{L^2}^2\right)^4 \|\na c_i\|_{L^2}^2 + \frac{D}{4} \|\Delta c_i\|_{L^2}^2
\end{aligned}
\ee using H\"older's inequality with exponents 6,3,2, the 3D Gagliardo-Nirenberg interpolation inequality applied to the mean-free vector $\na c_i$, Young's inequality for products, the Poincar\'e inequality applied to the mean-free function $\rho + \tilde{\rho}$, 3D elliptic estimates, and the relation between the density $\rho$ and the ionic concentrations.  As for the electromigration term, we have 
\be 
\beg{aligned}
&\left| Dz_i \int_{\TT^3} \na \cdot (c_i \na \Phi) \Delta c_i dx \right|
\\&\quad\quad\le D|z_i| \left(\|\na c_i\|_{L^2} \|\na \Phi\|_{L^{\infty}} + \|c_i - \bar{c}_i\|_{L^6}\|\Delta \Phi\|_{L^3} + |\bar{c}_i| \|\Delta \Phi\|_{L^2} \right)\|\Delta c_i\|_{L^2}
\\&\quad\quad\le \frac{D}{4}\|\Delta c_i\|_{L^2}^2 + C\|\na c_i\|_{L^2}^2 \|\na(\rho 
+ \tilde{\rho})\|_{L^2}^2  
\\&\quad\quad\le \frac{D}{4} \|\Delta c_i\|_{L^2}^2 + C\|\na c_i\|_{L^2}^2 \left(\|\na \tilde{\rho} \|_{L^2}^2 + \sum\limits_{j=1}^{N} \|\na c_i\|_{L^2}^2 \right)
\end{aligned}
\ee where the 3D elliptic estimate
\be 
\|\na \Phi\|_{L^{\infty}} + \|\Delta \Phi\|_{L^3} + \|\Delta \Phi\|_{L^2}
\le C\|\na \rho + \na \tilde{\rho} \|_{L^2}
\ee is exploited. Summing the equations \eqref{localco} over all indices $i \in \left\{1, \dots, N\right\}$, we infer that the energy inequality 
\be 
\frac{d}{dt} \sum\limits_{i=1}^{N} \|\na c_i\|_{L^2}^2 + D\sum\limits_{i=1}^{N} \|\Delta c_i\|_{L^2}^2 
\le C\left(\sum\limits_{i=1}^{N}\|\na c_i\|_{L^2}^2 \right)^4 + C_{\tilde{\rho}}
\ee holds, where $C_{\tilde{\rho}}$ is a positive constant depending only on $\|\na \tilde{\rho}\|_{L^2}$. This yields a local solution $(c_1, \dots, c_N)$ on a short time interval $[0,T_0]$. In order to prove uniqueness, suppose $(c_1^1, \dots, c_N^1)$ and $(c_1^2, \dots, c_N^2)$ are two solutions with equal initial data. For $i \in \left\{1, \dots, N \right\}$, let $c_i = c_i^1 - c_i^2$, $\rho = \rho_1 - \rho_2, \Phi = \Phi_1 - \Phi_2$, and $u = u_1 - u_2$, where $(\rho_1, \Phi_1, u_1)$ and $(\rho_2, \Phi_2, u_2)$ are the tuples of corresponding densities, potentials, and velocities associated with the two solutions. The potential difference $\Phi$ and velocity difference $u$
satisfy 
\be 
-\Delta \Phi = \rho =\sum\limits_{i=1}^{N} z_i c_i
\ee and 
\be 
u = -\PP(\rho \na \Phi_1 +(\rho_2 + \tilde{\rho})\na \Phi)
\ee respectively. 
Each $c_i$ evolves in time according to 
\be 
\pa_t c_i + D \Delta c_i 
= - u \cdot \na c_i^1 - u_2 \cdot \na c_i + Dz_i \na \cdot (c_i \na \Phi_1) + Dz_i \na \cdot (c_i^2 \na \Phi),
\ee which gives rise to the energy equation 
\be 
\frac{1}{2} \frac{d}{dt} \sum\limits_{i=1}^{N} \|c_i\|_{L^2}^2
+ D\sum\limits_{i=1}^{N} \|\na c_i\|_{L^2}^2
= - \sum\limits_{i=1}^{N} \int_{\TT^3} u \cdot \na c_i^1 c_i dx
+ D\sum\limits_{i=1}^{N} z_i \int_{\TT^3} (c_i \na \Phi_1 + c_i^2 \na \Phi ) \cdot \na c_i dx
\ee after multiplying by $c_i$, integrating over $\TT^3$, and summing over all indices $i = 1, \dots, N$. We bound the $L^2$ norm of $u$ as follows,
\be 
\beg{aligned}
\|u\|_{L^2}
&\le C\|\rho \na \Phi_1\|_{L^2} + C\|(\rho_2 +\tilde{\rho}) \na \Phi\|_{L^2}
\le C\|\rho\|_{L^2} \|\na \Phi_1\|_{L^{\infty}} + C\|\na \Phi\|_{L^6} \|\rho_2 + \tilde{\rho}\|_{L^3}
\\&\le C\|\rho\|_{L^2} \left(\|\rho_1 + \tilde{\rho}\|_{H^1} + \|\rho_2 + \tilde{\rho}\|_{H^1} \right)
\le C\left(\sum\limits_{i=1}^{N} \|c_i\|_{L^2} \right) \left(\|\rho_1 + \tilde{\rho}\|_{H^1} + \|\rho_2+ \tilde{\rho}\|_{H^1} \right)
\end{aligned}
\ee where the boundedness of $\PP$ in $L^2$, elliptic regularity, and Sobolev embeddings are used. Consequently, the nonlinear term in $u$ can be controlled by 
\be 
\beg{aligned}
\left|\int_{\TT^3} u \cdot \na c_i^1 c_i dx\right| 
&= \left|\int_{\TT^3} uc_i^1 \cdot \na c_i dx\right|
\le \|u\|_{L^2} \|c_i^1\|_{L^{\infty}}\|\na c_i\|_{L^2}
\\&\le \frac{D}{4}\|\na c_i\|_{L^2}^2
+ C\|c_i^1\|_{L^{\infty}}^2 \left(\|\rho_1\|_{H^1}^2 + \|\rho_2\|_{H^1}^2 + \|\tilde \rho \|_{H^1}^2 \right) \left(\sum\limits_{i=1}^{N} \|c_i\|_{L^2}^2 \right)
\\&\le \frac{D}{4}\|\na c_i\|_{L^2}^2
+ C\|c_i^1\|_{H^2}^2 \left(\|\rho_1\|_{H^1}^2 + \|\rho_2\|_{H^1}^2 + \|\tilde \rho \|_{H^1}^2 \right) \left(\sum\limits_{i=1}^{N} \|c_i\|_{L^2}^2 \right)
\end{aligned}
\ee after integrating by parts, using the divergence-free property of $u$, and applying the 3D continuous Sobolev embedding of $H^2$ in $L^{\infty}$. Elliptic estimates yield 
\be 
\beg{aligned}
&\left|D\sum\limits_{i=1}^{N} z_i \int_{\TT^3} (c_i \na \Phi_1 + c_i^2 \na \Phi ) \cdot \na c_i dx\right|
\le \frac{D}{4} \sum\limits_{i=1}^{N} \|\na c_i\|_{L^2}^2 + C\left(\|\na \Phi_1\|_{L^{\infty}}^2 + \|c_i^2\|_{L^3}^2 \right) \sum\limits_{i=1}^{N} \|c_i\|_{L^2}^2
\\&\quad\quad\quad\quad\le \frac{D}{4} \sum\limits_{i=1}^{N} \|\na c_i\|_{L^2}^2 + C\left(\sum\limits_{i=1}^{N} \|c_i^1\|_{H^1}^2 + \|c_i^2\|_{H^1}^2 + \|\tilde{\rho}\|_{H^1}^2 \right) \sum\limits_{i=1}^{N} \|c_i\|_{L^2}^2.
\end{aligned}
\ee Thus, it holds that 
\be \la{uncon}
\frac{d}{dt}\sum\limits_{i=1}^{N} \|c_i\|_{L^2}^2 
\le C\mathcal{G}(t) \sum\limits_{i=1}^{N} \|c_i\|_{L^2}^2
\ee where
\be 
\mathcal{G}(t) = \sum\limits_{i=1}^{N} \|c_i^1\|_{H^2}^2 \left(\|\rho_1\|_{H^1}^2 + \|\rho_2\|_{H^1}^2 + \|\tilde{\rho}\|_{H^1}^2 \right) + \sum\limits_{i=1}^{N} \|c_i^1\|_{H^1}^2 + \|c_i^2\|_{H^1}^2 + \|\tilde{\rho}\|_{H^1}^2 
\ee is integrable in time over $[0,T_0]$. By Gronwall's inequality, we deduce that 
\be 
\sum\limits_{i=1}^{N} \|c_i(t)\|_{L^2}^2
\le C\exp\left\{\mathcal{G}(t)\right\} \sum\limits_{i=1}^{N} \|c_i(0)\|_{L^2}^2
\ee for all $t \in [0,T_0]$. Since $c_i(0) = 0$, we deduce that $c_1(t) = c_2(t)$ for a.e. $x \in \TT^3$ and $t \ge 0$. This gives the uniqueness of solutions. 
\end{proof}

\section{Continuity of the Solution Map} \la{AP2}

In this appendix, we prove Proposition \ref{contsol}.

\beg{proof}
We denote the solutions of the NPD system corresponding to $w_0^1$ and $w_0^2$ respectively by $(c_1^1, \dots, c_N^1)$ and $(c_1^2, \dots, c_N^2)$, and we consider the following relations
\be 
\beg{aligned}
&\rho_1 = \sum\limits_{i=1}^{N} z_ic_{i}^{1},\rho_2 = \sum\limits_{i=1}^{N} z_ic_{i}^{2},  \Phi_1 = (-\Delta)^{-1}(\rho_1+\tilde{\rho}), \Phi_2 = (-\Delta)^{-1}(\rho_2+\tilde{\rho}), 
   \\&u_1 = - \mathbb{P} ((\rho_1+\tilde{\rho})\nabla \Phi_1),
    u_2 = - \mathbb{P} ((\rho_2+\tilde{\rho})\nabla \Phi_2).
\end{aligned}
\ee We consider the differences $c_i = c_{i}^{1} - c_{i}^{2}$ for $i\in \{1, \dots ,N\}, \rho = \rho_1 - \rho_2, \Phi = \Phi_1 - \Phi_2,$ and $u = u_1 - u_2.$ 
The concentration difference $c_i$ evolves according to \be\label{starequation}
\partial_tc_i + u_1 \cdot \nabla c_i + u \cdot \nabla c_i^{2} - D \Delta c_i = D z_i \nabla \cdot (c_{i}^{1} \nabla \Phi + c_i \nabla \Phi_2). 
\ee
Multiplying \eqref{starequation} by $-\Delta c_i$ and integrating over $\mathbb{T}^{3}$, we obtain 
\be
\beg{aligned}
    &\frac{1}{2} \frac{d}{dt} \| \nabla c_i \|_{L^2}^{2} + D \| \Delta c_i \|_{L^2}^{2} 
   \\&\quad\quad= \int_{\mathbb{T}^{3}} u_1 \cdot \nabla c_i\Delta c_i dx + \int_{\TT^3} u \cdot \nabla c_i^{2} \Delta c_i dx
    - \int_{\TT^3} D z_i \nabla \cdot (c_{i}^{1}\nabla \Phi + c_i\nabla \Phi_2)\Delta c_i dx 
    \\&
    \leq \|u_1\|_{L^{\infty}} \|\nabla c_i \|_{L^2}\|\Delta c_i\|_{L^2} + \|u \|_{L^3} \|\nabla c_i^{2} \|_{L^6}\|\Delta c_i\|_{L^2}
    \\&\quad\quad
    + D|z_i| \left[ \|\nabla c_{i}^{1} \|_{L^6} \|\nabla \Phi \|_{L^3} + \| c_{i}^{1} \|_{L^3} \|\Delta \Phi \|_{L^6} + \|\nabla c_{i} \|_{L^2} \|\nabla \Phi_2 \|_{L^{\infty}} + \| c_{i} \|_{L^2} \|\Delta \Phi_2 \|_{L^{\infty}} \right] \|\Delta c_i \|_{L^{2}}.   
\end{aligned}
\ee Since $H^2$ is continuously embedded in  $L^{\infty}$, the projector $\mathbb{P}$ is bounded in $H^2,$ and $H^2$ is a Banach Algebra, we have \be
    \|u_1 \|_{L^{\infty}} \leq C \|u_1 \|_{H^2}\leq C \|\rho_1 + \tilde{\rho}\|_{H^2} \|\nabla \Phi \|_{H^2} \leq C \|\rho_1 + \tilde{\rho}\|_{H^2} \|\nabla \rho \|_{L^2}.
\ee As $u = - \mathbb{P}(\rho \nabla \Phi_1 + (\rho_2 + \tilde{\rho})\nabla \Phi)$ and $\PP$ is bounded on $L^3$, it holds that \be
    \|u \|_{L^{3}} \le C \|\rho \|_{L^{6}}\|\na \Phi_1 \|_{L^{6}} + \| \rho_2 +\tilde{\rho} \|_{L^6} \|\nabla \Phi \|_{L^6}\\
    \leq C  \|\nabla \rho\|_{L^2} \left(\| \rho_1 + \tilde{\rho} \|_{L^2} + \| \rho_2 + \tilde{\rho} \|_{H^1}\right).
\ee By making use of Young's inequality for products, Sobolev estimates, and elliptic regularity, we deduce the differential inequality 
\be
\beg{aligned}
    &\frac{d}{dt} \|\nabla c_i\|_{L^2}^{2} + D \|\Delta c_i\|_{L^2}^{2}\\
    &\quad\quad\leq C \|\rho_1 + \tilde{\rho} \|_{H^2}^{2} \|\nabla \rho \|_{L^2}^{2} \|\nabla c_i \|_{L^2}^{2} + C\| \nabla \rho \|_{L^2}^{2} (\| \rho_1 + \tilde{\rho} \|_{L^2}^{2} + \| \rho_2 + \tilde{\rho} \|_{H^1}^{2})\|c_{i}^{2} \|_{H^2}^{2} \\
    &\quad\quad\quad\quad+ C\|c_i^1 \|_{H^2}^{2} \|\rho \|_{L^2}^{2} + C\|c_i^1 \|_{H^1}^{2} \|\nabla \rho \|_{L^2}^{2} +C \|\nabla c_i \|_{L^2}^{2} \| \rho_2 + \tilde{\rho} \|_{H^1}^{2} + C\| c_i \|_{L^2}^{2} \| \rho_2 + \tilde{\rho} \|_{H^2}^{2}.
    \end{aligned}
    \ee
    Summing over all indices $i \in \left\{1, \dots, N \right\}$, we infer that  \be 
    \label{equationnumber1}
    \beg{aligned}
       & \frac{d}{dt} \sum\limits_{i=1}^{N} \|\nabla c_i \|_{L^2}^{2} \\&\quad\quad\leq C\left(1 + \|\tilde{\rho} \|_{H^2}^{2} + \sum\limits_{i=1}^{N} \|c_i^1 \|_{H^2}^{2}+\sum\limits_{i=1}^{N} \|c_i^2 \|_{H^2}^{2}\right)\left(1 + \|\tilde{\rho} \|_{H^2}^{2} + \sum\limits_{i=1}^{N} \|c_i^1 \|_{H^1}^{2}+\sum\limits_{i=1}^{N} \|c_i^2 \|_{H^1}^{2}\right) \sum\limits_{i=1}^{N} \|c_i \|_{H^1}^{2}.
        \end{aligned}
    \ee As shown in \eqref{uncon}, it holds that 
    \be \la{equationnumber2}
    \frac{d}{dt}\sum\limits_{i=1}^{N} \|c_i\|_{L^2}^2 
\le C\mathcal{G}(t) \sum\limits_{i=1}^{N} \|c_i\|_{L^2}^2
    \ee Combining the energy inequalities  \ref{equationnumber1} and \ref{equationnumber2}, we deduce that \be
    \frac{d}{dt} \sum\limits_{i=1}^{N} \| c_i \|_{H^1}^{2} \leq CK(t) \sum\limits_{i=1}^{N} \|c_i \|_{H^1}^{2},
\ee  where \be
    K(t) =  \left(1 + \|\tilde{\rho} \|_{H^2}^{2} + \sum\limits_{i=1}^{N} \|c_i^1 \|_{H^2}^{2}+\sum\limits_{i=1}^{N} \|c_i^2 \|_{H^2}^{2}\right)\left(1 + \|\tilde{\rho} \|_{H^2}^{2} + \sum\limits_{i=1}^{N} \|c_i^1 \|_{H^1}^{2}+\sum\limits_{i=1}^{N} \|c_i^2 \|_{H^1}^{2}\right).   
\ee An application of Gronwall's inequality yields the desired instantaneous Lipschitz continuity in the norm of $\mathcal{V}$.

\end{proof}

\section{Decay of Volume Elements} \la{AP3}

In this appendix, we prove Proposition \ref{dve}.

\beg{proof}
We denote by $I$ the identity operator and  consider the operators $\mathcal{A}_{w_0}$ and $(\mathcal{A}_{w_0})_{n}$ defined by 
\beg{align*}
\mathcal{A}_{w_0}[\xi^j] &= (\mathcal{A}^1_{w_0}[\xi^j], \dots, \mathcal{A}^N_{w_0}[\xi^j]), \mathrm{where} \\\mathcal{A}^i_{w_0}[\xi^j]
&=  \PP(R_j \na \Phi + (\rho + \tilde{\rho})\na\Psi_j)  \cdot  \na c_i + u \cdot \na \xi_i^j - Dz_i \na \cdot (\xi_i^j \na \Phi + c_i \na \Psi_j),
\end{align*} and 
\be 
(\mathcal{A}_{w_0} - D\Delta)_{n} = (\mathcal{A}_{w_0} -D\Delta) \wedge I \wedge \dots \wedge I + \dots + I \wedge \dots \wedge I \wedge (\mathcal{A}_{w_0}-D\Delta)
\ee respectively. The evolution equation 
\be 
\pa_t (\xi_1 \wedge \dots \wedge \xi_n) + (\mathcal{A}_{w_0}-D\Delta)_n (\xi_1 \wedge \dots \wedge \xi_n) = 0,
\ee holds and yields 
\be 
\frac{d}{dt} V_n + Trace((\mathcal{A}_{w_0}-D\Delta)Q_n) V_n = 0
\ee where $Q_n$ is the orthogonal projection  in $\mathcal{V}$ onto the space spanned by $\xi^1, \dots, \xi^n$. Consequently, we infer that
\be 
V_n(t) \le V_N(0) \exp \left\{-\int_{0}^{t} Trace((\mathcal{A}_{w_0}-D\Delta) Q_n) ds \right\}
\ee for any $t \ge 0$. For each $t \ge 0$, we let $\left\{\phi_1, \dots, \phi_n \right\}$ be a orthonomal set in $\mathcal{V}$ spanning the linear span of $\left\{\xi^1, \dots, \xi^n \right\}$. Then 
\be 
\beg{aligned}
Trace((\mathcal{A}_{w_0}-D\Delta) Q_n) = \sum\limits_{j=1}^{n} (-D\Delta \phi_j + \mathcal{A}_{w_0}[\phi_j], -\Delta \phi_j)_{L^2}.
\end{aligned}
\ee Let $\mu_1, \dots, \mu_n$ be the first $n$ eigenvalues of the Laplace opeartor $-\Delta$ with periodic boundary conditions. On one hand, we have 
\be 
\sum\limits_{j=1}^{n} (-D\Delta \phi_j, -\Delta \phi_j)_{L^2}  \ge \mu_1 + \dots + \mu_n \ge CDn^2.
\ee On the other hand, we have 
\be 
\beg{aligned}
&\sum\limits_{j=1}^{n} |(\mathcal{A}_{w_0} [\phi_j], -\Delta \phi_j)_{L^2}|
\le \sum\limits_{j=1}^{n} \| \mathcal{A}_{w_0} [\phi_j]\|_{L^2} \|\Delta \phi_j\|_{L^2}
\\&\quad\quad\le C\sum\limits_{j=1}^{n} \sum\limits_{i=1}^{N} \|\PP(R_j \na \Phi + (\rho + \tilde{\rho}) \na \Psi_j\|_{L^3}\|\na c_i\|_{L^6} \|\Delta \phi_j\|_{L^2} + \|u\|_{L^{\infty}} \|\na \phi_i^j\|_{L^2} \|\Delta \phi_j\|_{L^2}
\\&\quad\quad\quad\quad+ CD\sum\limits_{j=1}^{n}\sum\limits_{i=1}^{N} \|\na \cdot (\phi_i^j \na \Phi + c_i \na \Psi_j) \|_{L^2} \|\Delta \phi_j\|_{L^2}
\end{aligned}
\ee where $\phi_j = (\phi_1^j, \dots, \phi_N^j), R_j = \sum\limits_{i=1}^{N} z_i \phi_i^j,$ and  $\Psi_j = (-\Delta)^{-1} R_j$. We estimate
\be 
\beg{aligned}
&\|\PP(R_j \na \Phi + (\rho + \tilde{\rho}) \na \Psi_j\|_{L^3}
\le C\left(\|R_j \na \Phi\|_{L^3} + \|(\rho + \tilde{\rho}) \na \Psi_j\|_{L^3} \right)
\\&\quad\quad\le C\|R_j\|_{L^6}\|\na \Phi\|_{L^6} + C\|\rho + \tilde{\rho}\|_{L^6} \|\na \Psi_j\|_{L^6}
\\&\quad\quad\le C\|R_j\|_{H^1} \|\rho + \tilde{\rho}\|_{L^2} + C\|\rho + \tilde{\rho}\|_{H^1} \|R_j\|_{L^2}
\\&\quad\quad\le C\|\rho + \tilde{\rho}\|_{H^1}
\end{aligned}
\ee using the Sobolev inequality, elliptic estimates, and the normalization $\|\phi^j\|_{\mathcal{V}} = 1$. Similarly, we have
\be 
\beg{aligned}
&\|\na \cdot (\phi_i^j \na \Phi + c_i \na \Psi_j)\|_{L^2}
\\&\quad\quad\le \|\na \phi_i^j\|_{L^2} \|\na \Phi\|_{L^{\infty}} + \|\phi_i^j\|_{L^6}\|\Delta \Phi\|_{L^3} + \|\na c_i\|_{L^2}\|\na \Psi_j\|_{L^{\infty}} + \|c_i\|_{L^3}\|\Delta \Psi_j\|_{L^6}
\\&\quad\quad\le C\|\rho + \tilde{\rho}\|_{H^1} + C\|c_i\|_{H^1}.
\end{aligned}
\ee Since $H^2$ is a Banach Algebra, it holds that 
\be 
\|u\|_{L^{\infty}}
\le C\|u\|_{H^2}
\le C\|(\rho +\tilde{\rho})\na \Phi\|_{H^2}
\le C\|\rho + \tilde{\rho}\|_{H^2} \|\na \Phi\|_{H^2}
\le C\|\rho + \tilde{\rho}\|_{H^2}\|\rho + \tilde{\rho}\|_{H^1}.
\ee Consequently, we infer that
\be 
\beg{aligned}
&\sum\limits_{j=1}^{n} |(\mathcal{A}_{w_0} [\phi_j], -\Delta \phi_j)_{L^2}|
\\&\quad\quad\le C\sum\limits_{j=1}^{n} \sum\limits_{i=1}^{N} \left(\|\rho + \tilde{\rho}\|_{H^1} \|\Delta c_i\|_{L^2} + \|\rho + \tilde{\rho}\|_{H^1} (\|\rho + \tilde{\rho}\|_{H^2}  +1) + \|c_i\|_{H^1} \right)\|\Delta \phi_j\|_{L^2}
\\&\quad\quad\le Cn^{\frac{1}{2}} \left(\sum\limits_{j=1}^{N}\|\Delta \phi_j\|_{L^2}^2\right)^{\frac{1}{2}} \left(\|\rho + \tilde{\rho}\|_{H^1} \left(\sum\limits_{i=1}^{N}\|\Delta c_i\|_{L^2} +\|\rho + \tilde{\rho}\|_{H^2}\right) + \sum\limits_{i=1}^{N} \|c_i\|_{H^1} + 1  \right).
\end{aligned}
\ee Since $(c_1, \dots, c_N) \in X$, there exists a constant $C(R)$ depending only on $R$ such that 
\be 
\left(\|\rho + \tilde{\rho}\|_{H^1} \left(\sum\limits_{i=1}^{N}\|\Delta c_i\|_{L^2} +\|\rho + \tilde{\rho}\|_{H^2}\right) + \sum\limits_{i=1}^{N} \|c_i\|_{H^1} + 1  \right)^2 \le C(R). 
\ee Thus,
\be 
\beg{aligned}
&\int_{0}^{t} Trace((\mathcal{A}_{w_0}-D\Delta) Q_n) ds
\\&\quad\ge CDn^2t - \frac{CDn^2}{2}t - Cn\int_{0}^{t} \left(\|\rho + \tilde{\rho}\|_{H^1} \left(\sum\limits_{i=1}^{N}\|\Delta c_i\|_{L^2} +\|\rho + \tilde{\rho}\|_{H^2}\right) + \sum\limits_{i=1}^{N} \|c_i\|_{H^1} + 1  \right)^2 ds
\\&\quad\ge CDn^2t - C(R)nt 
\ge CDn^2t
\end{aligned}
\ee provided that $n \ge C(R)$.
\end{proof}

\bibliographystyle{plain}
\bibliography{reference}

\end{document}